\documentclass[twocolumn]{autart}    

\usepackage{graphicx}          
\usepackage[cmintegrals]{newtxmath}
\usepackage{bm}
\usepackage{amsmath}
\usepackage{graphics}
\usepackage{color}
\usepackage{setspace}
\usepackage{enumerate}
\newtheorem{theorem}{Theorem}
\newtheorem{lemma}{Lemma}
\newtheorem{assumption}{Assumption}
\newtheorem{problem}{Problem}
\newtheorem{procedure}{Procedure}
\newtheorem{definition}{Definition}
\newtheorem{remark}{Remark}

\begin{document}

\begin{frontmatter}

\title{Power-Imbalance Allocation Control of Power Systems-Secondary Frequency Control} 

 \thanks[footnoteinfo]{Corresponding author: Kaihua Xi.}

\author{Kaihua Xi}\ead{K.Xi@tudelft.nl},    
\author{Johan L. A. Dubbeldam}\ead{J.L.A.Dubbeldam@tudelft.nl},              
\author{Hai Xiang Lin}\ead{H.X.Lin@tudelft.nl},
\author{Jan H. van Schuppen}\ead{jan.h.van.schuppen@xs4all.nl}

\address{Delft Institute of Applied Mathematics, Delft University of Technology, Mekelweg 4, 2628CD, Delft, The Netherlands}  

\begin{abstract}                          
The traditional secondary frequency control of power systems restores nominal frequency by steering \emph{Area 
Control Errors} (ACEs) to zero. Existing methods
are a form of integral control with the characteristic that 
large control gain coefficients introduce an overshoot and small ones result in a slow 
convergence to a steady state. In order to deal with the large frequency
deviation problem, which is the main concern of the power system integrated with a
large number of renewable energy, a faster convergence is critical. 
In this paper, we propose a secondary frequency control method
named \emph{Power-Imbalance Allocation Control} (PIAC) to restore the nominal frequency with a minimized control cost,
in which a coordinator estimates the power imbalance and dispatches the control inputs to the controllers
after solving an economic power dispatch problem. The power imbalance estimation converges exponentially in PIAC, both
overshoots and large frequency deviations are avoided. In addition, when PIAC is implemented in a multi-area controlled network,
the controllers of an area are independent of the disturbance of the neighbor areas, which allows an asynchronous control 
in the multi-area network. 
A Lyapunov stability analysis shows that PIAC is locally asymptotically stable and simulation results illustrates that it effectively
eliminates the drawback of the traditional integral control based methods. 
\end{abstract}

\end{frontmatter}

\section{Introduction}
%
%
%
%

Rapid expansion of the contribution of distributed renewable energy sources has
accelerated research efforts in controlling the power grid. In general, frequency control
is implemented at three different levels distinguished from fast to slow timescales \cite{Schavemaker2008,Ilic2000}. 
In a short time scale, the power grid is stabilized by decentralized droop control, which 
is called \emph{primary control}. While successfully balancing the power supply and demand, and synchronizing 
the power frequency, 
the primary control induces frequency deviations from the nominal frequency, e.g., 50 or 60 Hz.
The \emph{secondary frequency control} regulates the frequency back to its nominal frequency in a slower time scale 
than the primary control. On top of the primary and secondary control, the \emph{tertiary 
control} is concerned with global economic power dispatch over the networks in a large 
time scale. Consequently it depends on the energy prices and markets. 

The secondary frequency control is the focus of this paper. An interconnected electric power system 
can be described as a collection of subsystems, each of which is called a control area. 
The secondary control in a single area is regulated by \emph{Automatic Generation Control} (AGC),
which is driven by \emph {Area Control Error} (ACE). The ACE of an area is calculated from the 
local frequency deviations within the area and power transfers between the area and 
its neighbor areas. The AGC controls the power injections to force the ACE to zero, thus restores the nominal frequency.
 Due to the availability of a communication network, other secondary frequency control approaches have recently been developed 
 which minimize the control cost on-line \cite{Dorfler2014},  e.g., the Distributed
Average Integral Method (DAI) \cite{Zhao2015}, the Gather-Broadcast (GB) method \cite{Dorfler2016}, economic AGC (EAGC) method \cite{EAGC}, and
 distributed real time power optimal power control method \cite{Liu2016}.
These methods suffer from a common drawback, namely that they exhibit overshoot for large gain coefficients and
slow convergence for small gain coefficients \cite{overshoot_report2,overshoot_report1,Recent_survey_2005}. 
This is due to the fact that they rely on integral control which is well-known
to give rise to the two phenomena mentioned above. Note that the slow convergence speed results 
in a large frequency deviation which is the main concern of power systems integrated with a large amount of renewable energy. 

The presence of fluctuations is expected to increase in the near future, due to weather dependent 
renewable energy, such as solar and wind energy. These renewable 
power sources often cause serious frequency fluctuations and deviation from the nominal frequency due to the 
uncertainty of the weather.  This demonstrates 
the necessity of good secondary frequency control methods whose transient performance
is enhanced with respect to the traditional methods. We have recently derived such a method 
called \emph{Power Imbalance Allocation Method} (PIAC) in \cite{PIAC}, which can eliminate the 
drawback of the integral control based approach. This paper is the extended version of 
the conference paper \cite{PIAC} with additional stability analysis and the extension of PIAC to multi-area control.

We consider power systems with lossless transmission lines, which comprise
traditional synchronous machines, frequency dependent devices (e.g.,   
power inverters of renewable energy or freqeuncy dependent loads) and passive loads. 
We assume the system to be equipped with the primary controllers and propose the PIAC method in the framework of 
Proportional-Integral (PI) control, which first estimates the power imbalance of the system 
via the measured frequency deviations of the nodes of the synchronous machines and frequency 
dependent power sources, next dispatches the control inputs of the distributed controllers after solving the economic power 
dispatch problem. Since the estimated power imbalance 
converges exponentially at a rate that can be 
accelerated by increasing the control gain coefficient, the overshoot problem and
the large frequency deviation problem is avoided. Hence the drawback of 
the traditional ACE method is eliminated. Furthermore,
the control gain coefficient is independent to the 
parameters of the power system but only relies on the response time 
of the control devices. Consequently the transient performance is greatly 
enhanced by improving the performance of the control devices in PIAC. 
When implemented in a multi-area power network, PIAC makes the control 
actions of the areas independent, while the controllers of each area handles the power imbalance of the local area only.

 The paper is 
 organized as follows. We introduce the mathematical model of the power system in Section \ref{Sec:model}. We 
 formulate the problem and discuss the existing 
 approaches in Section \ref{Sec:control_frq1}, then propose the secondary frequency control approach, \emph{Power-Imbalance 
 Allocation Control} (PIAC), based on 
 estimated power imbalance in Section \ref{Sec:control_frq2} and analyze its the stability invoking the Lyapunov/LaSalle stability criterion 
 in Section \ref{Sec:stability}. Finally, we evaluate the performance of PIAC by simulations on the IEEE-39 New England test power 
 system in Section \ref{Sec:experiment}. Section \ref{Sec:conclusion} concludes with remarks. 

\section{The model}\label{Sec:model}


%
%
%

A power system is described by a graph $G=(\mathcal{V},\mathcal{E})$ with nodes $\mathcal{V}$ and edges
$\mathcal{E}\subseteq \mathcal{V}\times \mathcal{V}$, where a node represents a bus and edge $(i,j)$ represents
 the direct transmission line connection between node $i$ and $j$. 
 We consider a power system as a lossless 
electric network with constant voltage (e.g., transmission grids where the line resistances are neglected)
and an adjacency 
matrix $(\hat{B}_{ij})$
where $\hat{B}_{ij}$ denotes the susceptance between node $i$ and node $j$.
The system consists of three types of nodes, synchronous machines, frequency dependent devices  and 
passive loads, the sets of which are denoted by $\mathcal{V}_M$, $\mathcal{V}_F$ and $\mathcal{V}_P$ respectively. 
Thus $\mathcal{V}=\mathcal{V}_M\cup\mathcal{V}_P\cup\mathcal{V}_F$. The frequency dependent 
devices are for example frequency dependent loads, inverters of renewable energy, buses equipped with 
droop controllers.
Denote 
the number of the nodes in $\mathcal{V},\mathcal{V}_M,\mathcal{V}_F,\mathcal{V}_P$ by $n,~n_M,~n_F$, and $n_P$ respectively, hence 
$n=n_M+n_F+n_P$. The model is described by the following \emph{Differential Algebraic Equations} (DAEs), see e.g., \cite{Dorfler2016}, 
\begin{subequations}\label{eq:model1}
\begin{align}
\hspace{-10pt}
\dot{\theta}_i&=\omega_i,~i\in\mathcal{V}_M\cup\mathcal{V}_F,\\
M_{i}\dot{\omega}_{i}+D_{i}\omega_{i}&=P_{i}-\sum_{j\in \mathcal V}{B_{ij}\sin{(\theta_{i}-\theta_{j})}}+u_{i}, 
~i\in \mathcal{V}_M,\label{eq:syn}\\
D_{i}\omega_{i}&=P_{i}-\sum_{j\in \mathcal V}{B_{ij}\sin{(\theta_{i}-\theta_{j})}}+u_{i},~i\in \mathcal{V}_F,\label{eq:frq}\\
0&=P_{i}-\sum_{j\in \mathcal V}{B_{ij}\sin{(\theta_{i}-\theta_{j})}},~i\in \mathcal{V}_P, \label{eq:pass} 
\end{align}
\end{subequations}
where $\theta_i$ is the phase angle at node $i$, $\omega_i$ is the frequency 
deviation from the nominal frequency, i.e., $\omega_i=\overline{\omega}_i-f^*$ where 
$\overline{\omega}_i$ is the frequency and $f^*=50$ Hz or 60 Hz is the 
nominal frequency, $M_i>0$ denotes the moment of inertia of a synchronous machine,  $D_i>0$ is the droop control coefficient, 
$P_{i}$ is the power injection or demand,
$B_{ij}=\hat{B}_{ij}V_{i}V_{j}$ is the effective
susceptance of line $(i,j)$, $V_{i}$ is the voltage at node $i$, $u_i\in [\underline{u}_i,\overline{u}_i]$ is a secondary frequency control input.
Note that $u_i$ is a constrained input of the secondary frequency control, $\underline{u}_i$
and $\overline{u}_i$ are its lower and upper bound, respectively. Furthermore,
the set of nodes equipped with the secondary controllers is denoted by $ \mathcal{V}_K\subseteq \mathcal{V}_M\cup\mathcal{V}_F$ 
and $u_i=0$ for $i\notin  \mathcal{V}_K$. Here, we have assumed that the nodes that participate in secondary control
are equipped with primary controllers. Note that the loads can also be equipped with primary controllers \cite{Zhao2014}. 
The dynamics of the voltage and reactive power is not modeled, since they are 
irrelevant for control of the frequency. More details on decoupling the voltage and frequency control in the power system 
can be 
found in \cite{Kundur1994,Simpson-Porco2016,Trip2016240}.
The model with linearized sine functions in (\ref{eq:model1}) is also 
widely studied to design primary and secondary frequency control laws, 
e.g., \cite{Andreasson2013,EAGC,UC}. For the validity of the linearized model with lossless network, 
we refer to \cite{Ilic2000,Hertem2006}.

\section{Secondary frequency control of power systems}\label{Sec:control_frq1}

\subsection{Problem formulation}

In practice, the frequency deviation should be in a prescribed range in order to avoid damage to the devices in the power system. 
 We assume droop controllers to be installed at some nodes such that $\sum_{i\in \mathcal{V}_M\cup\mathcal{V}_F}{D_i}>0$. 
When the power supply and demand are time-invariant, the frequencies of all the nodes in $\mathcal{V}_M\cup\mathcal{V}_F$ synchronize
 at a state, called \emph{synchronous state} defined as follows,
 \begin{subequations}\label{eq:synchronize}
 \begin{align}
  \theta_i&=\omega_{syn}t+\theta_i^*,~i\in\mathcal{V},\\
  \omega_i&=\omega_{syn},~i\in\mathcal{V}_M\cup\mathcal{V}_F,\\
  \dot{\theta}_i&=\omega_{syn},~i\in\mathcal{V},\\
  \dot{\omega}_i&=0,~i\in\mathcal{V}_M\cup\mathcal{V}_F,
 \end{align}
 \end{subequations}
where $\omega_{syn}$ is the synchronized frequency deviation, and the phase angle 
differences at the steady state, $\{\theta_i^*-\theta_j^*,~(i,j)\in\mathcal{E}\}$, determine
the power flows in the transmission lines. 
The explicit synchronized frequency deviation $\omega_{syn}$ of the system is 
obtained by substituting (\ref{eq:synchronize}) into (\ref{eq:model1}) as 
\begin{eqnarray}
\hspace{50pt}\omega_{syn}=\frac{\sum_{i\in \mathcal{V}}{P_i}+\sum_{i\in\mathcal{V}_K}{u_i}}{\sum_{i\in \mathcal{V}_M \cup\mathcal{V}_F }{D_i}}. \label{eq:sync}
\end{eqnarray}
If and only if $\sum_{i\in \mathcal{V}}{P_i}+\sum_{i\in\mathcal{V}_K}{u_i}=0$, the frequency deviation of 
the steady state is zero, i.e., $\omega_{syn}=0$. 
 This implies that a system with only droop control, i.e., $u_i=0,\text{for}~ i\in\mathcal{V}_K$,  
can never converge to a steady state with 
$\omega_{syn}=0$ if the power demand and supply is unbalanced such that $\sum_{i\in \mathcal{V}}{P_i}\neq0$. 
This shows the need for the secondary control. As more renewable power sources are integrated 
into the power system, the power imbalance $\sum_{i\in \mathcal{V}}{P_i}$ may fluctuate severely leading to frequency oscillations. 
We aim to design an effective
method to control $\omega_{syn}$ to zero by solving the following problem, e.g., \cite{Dorfler2014} and \cite{Dorfler2016}. 

\begin{problem}\label{problem1}
Compute the inputs \{$u_i,i\in\mathcal{V}_K$\} of the power system so as to achieve the control objective of a 
balance of power supply and demand in terms of $\omega_{syn}=0$, or equivalently, $\sum_{i\in \mathcal{V}}{P_i}+\sum_{i\in\mathcal{V}_K}{u_i}=0$. 
\end{problem}

The following assumption states the basic condition for which a feasible solution of Problem \ref{problem1} exists.

\begin{assumption}\label{assumption1}
During a small time interval the value of power supply and demand are constant. In addition, for these values there exist control inputs
$\{u_i\in[\underline{u}_i,\overline{u}_i], i\in \mathcal{V}_K \}$, such that $\sum_{i\in\mathcal{V}}{P_i}+\sum_{i\in\mathcal{V}_K }{u_i}=0$. 
\end{assumption}

In a small time interval, the tertiary control, which calculates
the operating point stabilized by the secondary control, guarantees the existence of a steady state
and its local stability \cite{Ilic2000,Wood}. In addition, the main task of the generators
is to provide electricity to the loads and maintain the nominal frequency, so Assumption \ref{assumption1} is realistic.

For different controllers in the system, 
the control cost might be different for various reasons such as different device maintenance prices.
From the global perspective of the entire network, we aim to minimize the secondary frequency control cost which 
leads to Problem \ref{problem2}. 

\begin{problem}\label{problem2}
Compute the inputs \{$u_i,i\in\mathcal{V}$\} of the power system so as to achieve the control objective of 
minimal control cost, in addition to the control objective of a balance of power supply and demand with $\omega_{syn}=0$.  
\end{problem}

Corresponding to Problem \ref{problem2}, the following economic power dispatch problem needs to be solved, e.g.,\cite{Dorfler2016,Trip2016}.
\begin{eqnarray}
 \hspace{50pt}&&\min_{\{u_i,i\in\mathcal{V}_K\}} \sum_{i\in\mathcal{V}_K }J_i(u_i)
 \label{eq:optimal1}\\
 \nonumber
 &&s.t. ~~\sum_{i\in \mathcal{V}}{P_i}+\sum_{i\in\mathcal{V}_K }{u_i}=0. 
\end{eqnarray}
where  $J_i(u_i)$ is the control cost of node $i$, which incorporate the cost 
money and the constraints $u_i\in [\underline{u}_i,\overline{u}_i]$. 
Note that to solve (\ref{eq:optimal1}), the \emph{power imbalance} $\sum_{i\in \mathcal{V}}{P_i}$ should be 
known and the solution is for the small time interval mentioned in Assumption 1. 
Here, with respect to the existence of the solution of the economic power dispatch problem, we make the second assumption.

\begin{assumption}\label{assumption2}
The cost functions $J_i:\mathbb{R}\rightarrow \mathbb{R},~i\in\mathcal{V}_K $ are twice continuously differentiable and strictly convex 
such that $J''(u_i)>0$ where $J''(u_i)$ is the second order derivative of $J(u_i)$ with respect to $u_i$. 
\end{assumption}

Assumption \ref{assumption2} is also realistic because the constraint $u_i\in [\underline{u}_i,\overline{u}_i]$ 
can be incorporated in the objective function $J_i(u_i)$ for $i\in\mathcal{V}_K $ in a smooth way.

A necessary condition for solving the economic dispatch problem is \cite{Ilic2000},
\begin{eqnarray}
 \hspace{50pt} J'_i(u_i)=J'_j(u_j)=\lambda, ~ i,j,\in \mathcal{V}_K. 
 \label{eq:criterion}
\end{eqnarray}
where $J'_i(u_i)$ is the derivative 
of $J_i(u_i)$, which is the \emph{marginal cost} of node $i,i\in\mathcal{V}_K $, and $\lambda \in \mathbb{R}$ is the \emph{nodal price}. 
At the optimum of (\ref{eq:optimal1}) all the marginal costs 
of the controllers are equal to the nodal price. The \emph{market clearing price} $\lambda^*$ is obtained as the 
solution of the equation 
\begin{eqnarray}
\hspace{50pt} 0=\sum_{i\in\mathcal{V}}{P_i}+\sum_{i\in\mathcal{V}_K }{{J'}_i^{-1}(\lambda^*)}.
\end{eqnarray}
where ${J'}_i^{-1}(\cdot)$ is the inverse function of $J'_i(\cdot)$ which exists by Assumption \ref{assumption2}. 
Since in practice the power imbalance is uncertain with respect to the fluctuating power loads, 
the economic power dispatch problem (\ref{eq:optimal1}) cannot be solved directly.

\subsection{A brief review of secondary frequency control}

Before embarking on solving Problem \ref{problem1} and \ref{problem2}, we briefly outline existing
secondary frequency control methods and discuss their relevance for finding
a solution to Problem \ref{problem1} and \ref{problem2}.


\emph{ACE based AGC\cite{Ilic2000}:} The \emph{Area Control Error} (ACE) of an area 
is defined as 
\begin{eqnarray}
 \hspace{50pt} ACE=B\omega+P_{ex}-P_{ex}^*, \label{AGC_ACE}
\end{eqnarray}
where $B$ is a positive constant, $\omega$ is the frequency deviation of the area, $P_{ex}$ is the net power export, and $P_{ex}^*$ is 
the nominal value of $P_{ex}$. The adjustment of the power injection of the area is given 
as follows
\begin{eqnarray*}
\hspace{50pt} \dot{u}=-k\cdot ACE
\end{eqnarray*}
where $k\in (0,\infty)$ is a control gain coefficient. In the traditional \emph{Automatic Generation Control} (AGC) method, the frequency deviation is 
measured at a local node and communicated by a coordinator as the ACE to the controllers in the system, which calculate
their control inputs according to their participation factors. When the interconnected 
system is considered as a single area, the AGC has the form \cite{Dorfler2016} 
\begin{eqnarray}
 \hspace{50pt}\dot{\lambda}=-k\omega_{i^{*}},~ i\in\mathcal{V}, \quad u_i={J'}_i^{-1}(\lambda), \label{control:AGC}
\end{eqnarray}
where $k\in (0,+\infty)$ is a control gain coefficient, $\omega_{i^*}$ is the measured 
frequency deviation at a selected node $i^*$. $\lambda$ can be seen as the nodal price which 
converges to the market clearing price $\lambda^*$ as the power supply and demand is balanced. Note that the participation factor is involved 
in the derivative of the cost function, $J'_i(u_i)$. The frequency deviation of the area is not
well reflected in (\ref{control:AGC}) since it is 
measured at only one node. Furthermore, the communication network is not used so efficiently because
it only communicates the nodal price $\lambda$ from the coordinator to the controllers.

\emph{Gather-Broadcast (GB) Control\cite{Dorfler2016}:} In order to well reflect the frequency deviation of the area and 
use the communication network efficiently, the GB method measures the frequency deviations 
at all the nodes connected by the communication network. It has the form 
\begin{eqnarray}
 \hspace{40pt}\dot{\lambda}=-k\sum_{i\in\mathcal{V}}{C_i\omega_{i}},~ i\in\mathcal{V}, \quad u_i={J'}_i^{-1}(\lambda), 
 \label{eq:gather0}
\end{eqnarray}
where $k\in (0,+\infty)$ is a control gain coefficient and $C_i\in [0,1]$ is a set of convex weighting
coefficients with $\sum_{i\in\mathcal{V}}{C_i}=1$. The ACE in the GB method is actually the weighted average of the frequency deviations.  
As in the ACE based AGC method (\ref{control:AGC}),
a coordinator in the network also broadcasts the nodal price
to the controllers and the controllers compute the control inputs according to their own cost functions.

\emph{Distributed Averaging Integral control (DAI):} Unlike the ACE base AGC method and GB method implemented 
in a centralized way, DAI is implemented in a distributed way based on the consensus control principle \cite{Andreasson2013}. 
In the DAI method, there are no coordinators and 
each controller computes its own nodal price and communicates to its neighbors. A local node
in the system calculates its control input according to the local frequency deviation and the nodal prices received from its neighbors. 
As the state of the interconnected system reaches a 
new steady state, the nodal prices of all the nodes achieve a consensus at the market clearing price, thus 
Problem \ref{problem2} is solved. It has the form \cite{Zhao2015} 
\begin{eqnarray}
 \dot{\lambda}_i=k_i(-\omega_i+\sum_{j\in \mathcal{V}}w_{ij}(\lambda_j-\lambda_i)),\quad
 u_i={J'}_i^{-1}(\lambda_i), \label{eq:DAPI}
\end{eqnarray}
where $k_i\in (0,+\infty)$ is a gain coefficient for the controller $i$ and $\big(w_{ij}\big)$ denotes the undirected 
weighted communication network.
When $w_{ij}=0$ for all the lines of the communication network, the DAI control law reduces to 
a \emph{Decentralized Integral} (DecI) control law.
The DAI control has been widely studied on both 
the traditional power grids and Micro-Grids, e.g., \cite{Schiffer2016,Simpson-Porco2015}.
 Wu \emph{et.al.} \cite{XiangyuWu2016} proposed a distributed secondary control method where 
it is not necessary to know the nominal frequency for all the nodes in $\mathcal{V}_K$. 

When a steady state exists for the nonlinear system (\ref{eq:model1}), all the approaches above can restore the 
nominal frequency with an optimized control cost. However, it can be easily observed 
that the control approaches, e.g., (\ref{control:AGC}) and (\ref{eq:gather0}), are in the form of integral control where 
the control inputs are actually the integral of the frequency deviation. This 
will be further explained in Section \ref{Sec:control_frq2}. A common drawback of integral control is that 
\emph{the control input suffers from an overshoot problem with a large gain coefficient or  
slow convergence speed with a small one, which causes extra oscillation
or slow convergence of the frequency} \cite{overshoot_report2,overshoot_report1,Recent_survey_2005}.

The methods that we discussed above, concern controlling the nonlinear system (\ref{eq:model1}) based 
on passivity methods. 
However, the linearized version of the evolution equations (\ref{eq:model1}) was also addressed 
based on primal-dual method in the literature.
For example, Li \emph{et al.}  proposed an Economic AGC (EAGC) approach for the multi-area
frequency control \cite{EAGC}. In this method each controller exchanges control signals that are used
to successfully steer the state of the system to a steady state with optimized dispatch, by a partial
primal-dual gradient algorithm \cite{Boyd2010}.
Unfortunately, transient performance was not considered in the method. 
A potentially very promising method to the control of the linear system was recently 
proposed by Zhao \emph{et al.} In \cite{UC} a novel framework for primary and secondary control,
which is called Unified Control (UC), was developed. The advantage of UC is that it automatically
takes care of the congestions that may occur in the transmission lines. Numerical simulations
results show that the UC can effectively reduce the harmful large transient oscillations in the frequency.
However, so far a theoretical analysis on how the UC improves the transient performance is lacking.  
Other recently  reported study are by Liu \emph{et al.} in \cite{Liu2016} and 
by Z. Wang \emph{et al.} in \cite{Wangzhaojian2017}. These methods optimize both control costs and manages power flow congestion using the
principle of consensus control, but cannot prevent large frequency deviations at short times.
For more details of distributed frequency control, we refer to 
the survey paper \cite{Molzahn2017}.


Finally, we mention here control methods whose underlying principle is neither based on integral
control nor on primal-dual method. The optimal load-frequency control framework by Liu \emph{et al.}, described in \cite{Liu2003}, is an
example of such a method. The goal is still to optimize the control costs and  frequency deviation,
but rephrasing it as a finite horizon optimization problem. This can only be solved when the
loads are precisely known within the selected time horizon. Obviously, this will require very precise forecasting of the loads.
A more robust approach based on the concept of the Active Disturbance Rejection Control (ADRC) \cite{Han2009}, was
pursued by Dong \emph{et al.} \cite{Dong2012}.  The method is robust against model uncertainties, parameter variations
and large perturbations which was employed to construct a decentralized load frequency approach for interconnected systems. However,
the decentralized control employed in this method prevents a solution to Problem \ref{problem2}.

%
%

As more renewable power sources are integrated into the power system, the fluctuations in the power supply become faster and larger. 
There is a need to design a control law that has a good transient performance without the overshoot problem 
and with a fast convergence speed. The traditional method to eliminate the overshoot is to calculate the 
control gain coefficients by analyzing the eigenvalue of the linearized systems \cite{XiangyuWu2017,Wu2016}. However, 
the improvement of the transient performance obtained by the eigenvalue analysis is still poor because of 
the dependence of the eigenvalues on the control law structure, and the large scale, 
complex topology and heterogeneous power generation and loads of the power system.

Based on the framework of PI control this paper aims to design a secondary frequency control method that
remedies the drawbacks mentioned above of the existing methods. To this end, we consider 
the following problem concerning the 
transient performance of the power system (\ref{eq:model1}) after a disturbance.
\begin{problem}\label{problem3}
 For the power system (\ref{eq:model1}) with primary and secondary controllers, 
 design a secondary frequency control law for 
 $\{u_i,~i\in\mathcal{V}_K\}$ so as to eliminate the extra oscillation
 of the frequency caused by the overshoots of the control inputs, thus improve
 the transient performance of the system after a disturbance through
 an accelerated convergence of the control inputs. 
\end{problem}

To address Problem \ref{problem3}, various control laws have been proposed for the frequency control of power systems, e.g., 
the sliding mode based control laws \cite{overshoot2,overshoot3} and $\mathcal{H}_2/\mathcal{H}_{\infty}$ control based control laws \cite[Chapter 5]{Hassan_Bevrani} and \cite{Hinfinty_control_overshoot} \emph{et al},
which are able to shorten the transient phase while avoiding the overshoots and large frequency deviations.   
However, they focus on the linearized system and did not consider the economic power dispatch problem (\ref{eq:optimal1}) at the steady state.

%

\section{Power Imbalance Allocation Control}\label{Sec:control_frq2}

In this section, we introcuce the \emph{ Power-Imbalance Allocation Control} (PIAC) method
to solve Problem \ref{problem1}-\ref{problem3}. 

The communication network is necessary for solving Problem 2, 
for which we make the following assumption. 

\begin{assumption}\label{assumption3} All the 
buses in $\mathcal{V}_M\cup \mathcal{V}_F$ 
can communicate with a coordinator at a central location via a communication network. 
The frequency deviations $\omega_i,~i\in \mathcal{V}_M\cup \mathcal{V}_F$, can 
be measured and subsequently communicated to the coordinator. For ${i\in\mathcal{V}_K}$,
the control input $u_i$ can be computed by the coordinator and dispatched to
the controller at node $i$ via the communication network. 
\end{assumption}

In Assumption \ref{assumption3}, the local nodes need to provide the coordinator the local frequency deviation which 
are the differences between the measured frequencies by meters and the nominal frequency. We remark that there 
usually are time-delays in the measurement and communications which are neglected in this paper.

In the following, we first define an abstract frequency deviation of the power system, for which 
we design a control law to eliminate the overshoot of the control input, then introduce the 
 PIAC method firstly control  the power system controlled as a single area and followed by a  
multi-area approach in subsection (\ref{subsec:single}) and (\ref{subsec:multi_description}) respectively. 

\begin{definition}\label{definition1}
For the power system (\ref{eq:model1}), define an abstract frequency deviation $\omega_s$ such that 
 \begin{eqnarray}\label{eq:global}
 ~~~~~~~~~~~~~~M_s\dot{\omega}_s=P_s-D_s\omega_s+u_s,
\end{eqnarray}
where $M_s=\sum_{i\in\mathcal{V}_M}{M_i}$, $P_s=\sum_{i\in\mathcal{V}}{P_i}$, $D_s=\sum_{i\in\mathcal{V}}{D_i}$,
and $u_s=\sum_{i\in\mathcal{V}_K}u_i$.
\end{definition}

Note that $\omega_s(t)$ is a \emph{virtual global frequency deviation}, which neither equals 
to $\omega_i(t)$ nor $\omega_{syn}(t)$ in general. However, 
if the power loads (or generation) $P_i$ change 
slowly and the frequencies $\{\omega_i,~i\in\mathcal{V}_M\cup\mathcal{V}_F\}$ synchronize quickly,
the differences between $\omega_i(t)$ and $\omega_s(t)$ are 
negligible, i.e., $\omega_i(t)=\omega_s(t)$, by summing all the equations of (\ref{eq:model1}), 
(\ref{eq:global}) is derived.

At the steady state of (\ref{eq:global}), we have $\omega_s^*=(P_s+u_s)/D_s$ which leads to $\omega_s^*=\omega_{syn}^*$ by 
(\ref{eq:sync}). Then the objective of the secondary 
frequency control that restores the nominal frequency of the system (\ref{eq:model1}) is equivalent
to controlling $\omega_s$ of (\ref{eq:global}) to zero, i.e., 
$\lim_{t\rightarrow\infty}{\omega_s(t)=0}.$

Since the frequencies $\{\omega_i,~i\in\mathcal{V}_M\cup\mathcal{V}_F\}$ of 
system (\ref{eq:model1}) with primary controllers synchronize quickly, 
the extra oscillations of $\omega_i$ are actually 
introduced by the overshoot of the total amount of control inputs $u_s$ in the traditional control 
approach (\ref{control:AGC}). This 
is because it is in the form of integral control. It will be 
explained in an extreme case where system (\ref{eq:model1}) is well controlled
by primary controllers such that 
$\omega_i=\omega_s$ for all $i\in\mathcal{V}$. It
can be obtained easily for (\ref{control:AGC}) by substituting 
$\omega_i$ by $\omega_s$ that the total amount of control inputs, $u_s$, 
is calculated as follows,
\begin{subequations}\label{integral}
 \begin{align}
 &\dot{\lambda}_s=\omega_s,\\
 &u_s=-k\lambda_s,
 \end{align}
\end{subequations}
which is in the form of integral control. 
Note that (\ref{integral}) can 
also be obtained similarly for the GB method (\ref{eq:gather0}) 
and for the DAI method (\ref{eq:DAPI}) with a special setting of
control gain coefficients $k_i$ and communication weight $w_{ij}$ e.g., $\{k_i,~i\in\mathcal{V}_k\}$ are all identical
and $w_{ij}$ forms a Laplacian matrix, $L=(w_{ij})$, such that $\sum_{i=1}^{n_K}{w_{ij}}=0$ for all $j=1,\cdots,n_K$(
$n_K$ is the number of nodes in $\mathcal{V}_K$).

In order to accelerate the convergence of the frequency to its nominal value without introducing extra oscillations, 
the overshoot of $u_s$ should be avoided. 
Similar to the PI cruise control of a car \cite{PIDcontrol}, for $\omega_s$ in (\ref{eq:global}), 
we introduce the following control law, 
\begin{subequations}\label{eq:controlPIAC2}
 \begin{align}
   &\dot{\eta}_s=D_s\omega_s,\label{eq:controlPIAC2a}\\
 &u_s=-kM_s\omega_s-k\eta_s,
 \end{align}
\end{subequations}
where $k\in (0,+\infty)$ is a parameter to be chosen by engineers. Details on the setting of this 
parameter will be further discussed after introducing the PIAC method in this section. 
From (\ref{eq:controlPIAC2}) and (\ref{eq:global}), we obtain
\begin{eqnarray}
 ~~~~~~~~~~~~~~~~~~~~~\dot{u}_s(t)=-k(P_s+u_s(t)), \label{eq:piac_u}
\end{eqnarray}
which indicates that $u_s$ is an estimate of $-P_s$ and converges to the 
power imbalance $-P_s$ exponentially. Hence, the overshoot of $u_s$ is 
eliminated. With the value of $u_s$ obtained from (\ref{eq:controlPIAC2}), 
the control input $u_i$ for node $i\in\mathcal{V}_K$ 
is computed by solving the following economic power dispatch problem,  
\begin{eqnarray}
 ~~~~~~~~~~~~~&&\min_{\{u_i\in \mathbb{R},i\in\mathcal{V}_K\}} \sum_{i\in\mathcal{V}_K }J_i(u_i),
 \label{eq:optimal2}\\
 \nonumber
 &&s.t. ~~-u_s(t)+\sum_{i\in \mathcal{V}_K }{u_i}(t)=0. 
\end{eqnarray}
However, $u_s$ cannot be 
calculated as in (\ref{eq:controlPIAC2}) since $\omega_s$ is a 
virtual frequency deviation which cannot be measured in practice. 
In the following, we introduce the PIAC method where $u_s$ also converges to $-P_s$ exponentially as in (\ref{eq:piac_u}). 

%
%
%

\subsection{Single-area implementation of PIAC}\label{subsec:single}

We consider the power system controlled as a single area without 
any power export (or import). The PIAC method is defined 
as follows.
\begin{definition}[PIAC]\label{definition2}
Consider the power system described by (\ref{eq:model1}) with 
assumption (\ref{assumption1}-\ref{assumption3}), the PIAC control law is defined as the dynamic controller,
\begin{subequations}\label{eq:PIAC}
 \begin{align}
\dot{\eta}(t)&=\sum_{i\in\mathcal{V}_M\cup\mathcal{V}_F}D_i\omega_i(t), \label{eq:PIACa}\\
u_s(t)&=-k\Big(\eta(t)+\sum_{j\in\mathcal{V}_M}{M_j\omega_j(t)}\Big),\label{eq:PIACb}\\
0&=-u_s(t)+\sum_{i\in\mathcal{V}_K}{{J'}_i^{-1}(\lambda(t))}, \label{eq:PIACc}\\
u_i(t)&={J'}_i^{-1}(\lambda(t)),~i\in\mathcal{V}_K,  \label{eq:PIACd}
 \end{align}
\end{subequations}
where $k\in(0,\infty)$ is a parameter of the control law, $\eta$ is a state variable 
introduced for the integral term as $\eta_s$ in (\ref{eq:controlPIAC2a}), and $\lambda$ is an algebraic variable for solving the 
optimization problem (\ref{eq:optimal2}), which is a function of time. 
\end{definition}

For the special case with quadratic cost function $J_i(u_i)=\frac{1}{2}\alpha_iu_i^2$ for node $i$, 
the control law becomes
\begin{subequations}
 \begin{align*}
  \dot{\eta}&=\sum_{i\in\mathcal{V}_M\cup\mathcal{V}_F}{D_i\omega_i},\\
  u_s&=-k(\sum_{i\in\mathcal{V}_M}{M_i\omega_i+\eta_s}),\\
  u_i&=\frac{\alpha_s}{\alpha_i}u_s,~i\in\mathcal{V}_K,\\
  \frac{1}{\alpha_s}&=\sum_{i\in\mathcal{V}_K}\frac{1}{\alpha_i},
 \end{align*}
\end{subequations}
where $\alpha_i\in\mathbb{R}$ is the control price of node $i$. 

PIAC is based on the design principle of coordination. The local nodes in $\mathcal{V}_M\cup \mathcal{V}_F$
send the measured frequency deviations to the coordinator. The coordinator computes the control inputs
$\{u_i, i\in\mathcal{V}_K\}$ or marginal cost $\lambda$
using the local measurements and sends them to all local controllers of the nodes indexed by $ \mathcal{V}_K$. 
 The procedure is similar to the GB method which gathers the locally measured frequency deviation from the nodes and broadcasts the 
nodal price to the controllers.
The control procedure and properties of PIAC are summarized in the following procedure and theorem
respectively. 
\begin{procedure}\label{procedure}
 If the assumptions \ref{assumption1}-\ref{assumption3} hold, then the secondary frequency control approach of the power system (\ref{eq:model1}), PIAC, 
is implemented as follows.  
\begin{enumerate}[(i)]
 \item Collect the measured frequency deviations $$\{\omega_i,~i\in \mathcal{V}_M \cup\mathcal{V}_F \},$$
 \item Calculate the total amount $u_s(t)$ of control inputs by (\ref{eq:PIACb}),
 \item Solve the optimization problem (\ref{eq:optimal2}) by (\ref{eq:PIACc},\ref{eq:PIACd}), 
 \item Allocate the power compensation $\{u_i,~i\in \mathcal{V}_K \}$ to the controllers. 
\end{enumerate}
\end{procedure}

\begin{theorem}\label{theorem1}
 Consider the power system (\ref{eq:model1}) and the PIAC control law, the controller has
 the following properties,
 \begin{enumerate}[(a)]
\item at any time, $t\in T$, $u_s(t)$ satisfies (\ref{eq:piac_u}).
Thus it is an estimate of the power imbalance $-P_s$;
\item at any time, $t\in T$, the input values 
$\{u_i(t),i\in\mathcal{V}_K\}$ are computed by solving the
optimization problem (\ref{eq:optimal2}). 
So the total amount of control inputs, $u_s(t)$, are dispatched to the local nodes of $\mathcal{V}_K$ economically; 
\item at the steady state, $\omega_{i}^*=0$ for all $i\in\mathcal{V}_M\cup\mathcal{V}_F$ and the power-imbalance 
$P_s$ is optimally compensated by the local controllers. Hence both Problem \ref{problem1} and Problem \ref{problem2} 
are solved;
\item because 
  $$u_s(t)=-k\sum_{i\in\mathcal{V}_M}{M_i\omega_i}(t)-k\int_{0}^{t}{\sum_{j\in\mathcal{V}_M\cup\mathcal{V}_F}D_j\omega_j(\tau) d\tau},$$
the PIAC control law is of proportional-integral type. 
 \end{enumerate}
 \end{theorem}
\emph{Proof:} (a). From (\ref{eq:model1}) it follows,  
\begin{eqnarray}
\sum_{i\in\mathcal{V}_M}M_i\dot{\omega}_i=\sum_{i\in\mathcal{V}}P_i-\sum_{i\in\mathcal{V}_M\cup\mathcal{V}_F}{D_i\omega_i}+\sum_{i\in\mathcal{V}}{u_i}.
\label{eq:zt}
\end{eqnarray}
It follows from (\ref{eq:PIACa},\ref{eq:PIACb}) that,
\begin{subequations}
\begin{align*}
\dot{u}_s(t)&=-k\Big(\sum_{i\in\mathcal{V}_M\cup\mathcal{V}_F}D_i\omega_i(t)+\sum_{j\in\mathcal{V}_M}M_j\dot{\omega}_i(t)\Big)\\
&~\text{by (\ref{eq:zt}),}\\
&=-k\Big(\sum_{i\in\mathcal{V}_K}{u_i(t)}+\sum_{j\in\mathcal{V}}{P_i}\Big)\\
 &~\text{by (\ref{eq:PIACd}),}\\
&=-k\Big(\sum_{i\in\mathcal{V}_K}{{J'}_{i}^{-1}(\lambda)}+\sum_{j\in\mathcal{V}}{P_i}\Big)\\
 &~\text{by (\ref{eq:PIACc}),}\\
&=-k\Big(u_s(t)+P_s\Big),
\end{align*}
\end{subequations}
Thus (\ref{eq:piac_u}) is obtained which indicates that $u_s(t)$ is an estimate of $-P_s$ with $k\in(0,\infty)$, i.e.,
$u_s(t)$ converges to $-P_s$ exponentially with a speed determined by the control gain coefficient $k$.

(b). According to the definition of the PIAC control law that at any time, $t\in T$, 
 \begin{eqnarray}
 \nonumber
  ~~~~~~~~~~J'_i(u_i(t))=J'_j(u_j(t))=\lambda(t),~\forall~ i,j\in\mathcal{V}_K.
 \end{eqnarray}
Thus the necessary condition (\ref{eq:criterion}) for economic dispatch of $u_s(t)$ to all the local nodes is satisfied and 
at any time $t\in T$, the control inputs
$\{u_i,~i\in\mathcal{V}_K\}$ solve the optimization problem (\ref{eq:optimal2}). 
Hence, $\{u_i(t),i\in\mathcal{V}\}$ achieve the minimal cost. 

(c). It follows from (\ref{eq:piac_u}) that $u^*_s=-\sum_{i\in\mathcal{V}}P_i$ at the steady state, which 
yields $\omega_{syn}=0$ by (\ref{eq:sync}). Thus $\omega^*_i=0$ for all $i\in\mathcal{V}_M\cup\mathcal{V}_F$ by 
(\ref{eq:synchronize}), Problem \ref{problem1} is solved. 
It follows from (b) that the control inputs $\{u_i,i\in\mathcal{V}\}$ achieve
minimal cost at the steady state. So the optimization problem (\ref{eq:optimal1}) is solved
at the steady state. 
Hence, Problem \ref{problem2} is solved. 

(d). This follows directly from the definition of the PIAC control law. \hfill $\Box$

Note that in Theorem \ref{theorem1}(c), we have assumed that the steady state exists which will be
further described in Assumption \ref{assumption4} in Section 5. 

In order to clearly illustrate how PIAC improves the transient performance of the system,
with the abstract frequency deviation $\omega_s$ and control inputs $u_s$, we decompose
the dynamic process into three sub-processes,
\begin{enumerate}[(i)]
 \item the convergence process of $u_s$ to $-P_s$ as in (\ref{eq:piac_u}) with a speed 
 determined by $k$.
 \item the convergence process of $\omega_s$ to zero as in (\ref{eq:global}) with a speed determined by $u_s$ and $D_s$.  
 \item and the synchronization process
of $\{\omega_i,~i\in\mathcal{V}_M\cup\mathcal{V}_F\}$ to $\omega_s$ described by (\ref{eq:model1}) and 
with the synchronization speed mainly determined by $\{D_i,~i\in\mathcal{V}_M\cup\mathcal{V}_F\}$.
\end{enumerate}
The transient performance of the power system can be improved 
by tuning the corresponding 
control gain coefficients for the sub-processes. PIAC focuses on the first two sub-processes only in which the 
transient behaviors of $u_s$ and $\omega_s$ can be improved by increasing the parameter $k$, and 
the primary control focuses on the synchronization of $\{\omega_i,~i\in\mathcal{V}_M\cup\mathcal{V}_F\}$
which can be improved by tuning the coefficient $\{D_i,~i\in\mathcal{V}_M\cup\mathcal{V}_F\}$ 
as in \cite{Dorfler2014,Motter2013}. 
As $u_s$ converges to $-P_s$ exponentially, the control inputs $\{u_i,~i\in\mathcal{V}_K\}$ in (\ref{eq:PIACd})
also converge to the desired optimal values in an exponential way without any overshoots and their convergence
speed increase with a large $k$. Thus the extra oscillation caused by 
the overshoot of $u_s$ is avoided and Problem \ref{problem3} is solved by PIAC.

With the process decomposition, the improvement in the transient performance by the UC control law \cite{UC}
can be explained, and the control gain coefficients of the primal-dual based control laws, e.g.,\cite{Wangzhaojian2017}
without triggered constraints of line flows, which also estimate the power disturbance,
can be tuned to improve the transient performance.

In the following, we introduce the properties of PIAC in the
implementation for the power system controlled as a single area. 

On the communications in PIAC, note that the communicated data 
include dynamic data and static data. The dynamic data are 
the frequency deviations $\{\omega_i,~i\in\mathcal{V}_M\cup\mathcal{V}_F\}$ and the 
control inputs $\{u_i~i\in\mathcal{V}_K\}$, both should be communicated frequently. The static data are 
the moment inertia $\{M_i,~i\in\mathcal{V}_i\}$, the droop control coefficients $\{D_i,~i\in\mathcal{V}_M\cup\mathcal{V}_F\}$
and the cost functions $\{J_i(u_i),i\in\mathcal{V}_K\}$, which are constant 
in a short time as in Assumption \ref{assumption1} and are not necessary to communicate as frequently as the dynamic data. 

In the computing of the control inputs $\{u_i,~i\in\mathcal{V}_K\}$, 
solving the optimization problem (\ref{eq:optimal2}) is equivalent to solving the equations (\ref{eq:PIACc}) and (\ref{eq:PIACd}).
So the computation includes calculating the integral of frequency deviation in  (\ref{eq:PIACa}) and 
solving the algebraic equation (\ref{eq:PIACc},\ref{eq:PIACd}). For quadratic cost functions $J_i(u_i)=\frac{1}{2}\alpha_iu_i^2$, the 
computation requires approximately $4n_M+n_F$ arithmetic operations that are $\times$, +, - etc. For nonlinear cost functions, 
an iteration method is used to solve the one dimension nonlinear algebraic equation (\ref{eq:PIACc}) for $\lambda$, which 
needs more computing than for the quadratic cost functions. 
\begin{remark}
PIAC is a centralized control where the 
 communications and computations increase as the scale of the power systems increases. In 
 a large-scale power system, a large number of devices communicate to the coordinator simultaneously. In that case,
 the time-delay due to the communications and computing is not negligible, in such a situation
further investgation for the transient performance and stability of the nonlinear system is needed. 
\end{remark}

On the dynamics of $u_s$, it can be observed from (\ref{eq:piac_u}) that $u_s$ converges exponentially with a speed that 
can be accelerated by increasing the control gain coefficient $k$ which does 
not depend on any parameters of the power system and the cost functions. Hence, when the dynamics of 
the voltages are considered which does not 
influence the dynamics of $u_s$ in (\ref{eq:piac_u}), the power supply and demand can also be balanced. 

On the control gain coefficient $k$, we remark that it does neither 
depend on the parameters of the system nor on the economic power dispatch problem,
and can be set very large from the perspective of theoretical analysis. However,
it relies on how sensitive the control devices are to the fluctuation of 
power imbalance. In this case, it can be tuned according to the response time of the control devices and the desired range of the 
frequency deviation. The control actuators of traditional hydraulic and steam turbines are their 
governor systems. For details of the model of the governor system and its response time, we refer to \cite[chapter 9,11]{Kundur1994}.

By Assumption \ref{assumption3}, 
the PIAC method requires that all the nodes in $\mathcal{V}_M\cup\mathcal{V}_F$ can communicate with the 
 coordinator. However, our initial numerical experiments show that 
 PIAC can still restore the nominal frequency even without all the 
 frequency deviation collected, where the estimated power imbalance converges to the actual
 power imbalance although not exponentially. This is because PIAC includes the integral control which drives
 the synchronized frequency deviation to zero.

\begin{remark}
 In practice, the state of the power system is never at a true equilibrium state, 
because of the fluctuating of the power imbalance caused by the power loads. 
Furthermore, the fluctuations becomes even more serious when a large amount of renewable power sources 
are integrated into the power system.  
 In this case, it is more practical to model the power imbalance as a time-varying function. 
 For the power system with time-varying power imbalance, 
 analysis and numerical simulations show that PIAC 
 is also able to effectively control the synchronized frequency to any desired range by increasing the control gain coefficient $k$ \cite{PIAC2}. 
\end{remark}


\subsection {Multi-area implementation of PIAC}\label{subsec:multi_description}

When a power system is partitioned into several control areas, each area has either an export or an import of power. 
After a disturbance occurred, 
the power export or the power import of an area must be restored to the nominal value calculated by tertiary control. 
In this subsection, we introduce the implementation of PIAC for multi-area control where 
each area has power export (or import).

Denote the set of the nodes in the area $A_r$ by $\mathcal{V}_{A_r}$, the set of the boundary lines between area $A_r$ and all the 
other areas by $\mathcal{E}_{A_r}$.  Denote $\mathcal{V}_{M_r}=\mathcal{V}_M\cap\mathcal{V}_{A_r}$,
$\mathcal{V}_{F_r}=\mathcal{V}_F\cap\mathcal{V}_{A_r}$ and $\mathcal{V}_{K_r}=\mathcal{V}_K\cap\mathcal{V}_{A_r}$. The multi-area control of PIAC in the area $A_r$ is defined as 
follows. 

\begin{definition}
Consider the power system (\ref{eq:model1}) controlled as
several areas with Assumption (\ref{assumption1}-\ref{assumption3}), the multi-area implementation 
of PIAC for area $A_r$ is defined as the dynamic controllers,
\begin{subequations}\label{eq:mulit}
 \begin{align}
\dot{\eta}_r(t)&=\sum_{i\in \mathcal{V}_{M_r} \cup \mathcal{V}_{F_r} }D_i\omega_i(t)+P_{ex}(t)-P_{ex}^*,\label{eq:multi}\\
u_r(t)&=-k_r(\sum_{i\in \mathcal{V}_{M_r} }M_i\omega_i(t)+\eta_r(t)),\\
0&=-u_r+\sum_{i\in\mathcal{V}_{K_r}}{{J'}_i^{-1}(\lambda_r)},\label{eq:constraints_multi}\\
u_i&={J'}_i^{-1}{(\lambda_r)}, ~i\in\mathcal{V}_{M_r}\cup\mathcal{V}_{F_r}. 
 \end{align}
\end{subequations}
where $P_{ex}=\sum_{i\in\mathcal{V}_{A_r}, (i,j)\in \mathcal{E}_{A_r}}B_{ij}\sin{(\theta_{i}-\theta_{j})}$ is the export power of area $A_{r}$,  
$P_{ex}^*$ is the nominal value of $P_{ex}$, $k_r\in (0,\infty)$ is a control gain coefficient, $\eta_r$ 
is a state variable for the integral term as $\eta_s$ in (\ref{eq:controlPIAC2a}) and $\lambda_r$ is an algebraic variable for the 
controller. 
\end{definition}

It can be observed from (\ref{eq:mulit}) that the control procedure for the 
coordinator in area $A_r$ is similar to Procedure \ref{procedure} but the power export 
deviation, $P_{ex}-P_{ex}^*$, should be measured.
The sum of the three terms in the right hand side of (\ref{eq:multi})
is actually the ACE of the area. 
PIAC has the proportional control input
included in secondary frequency control, which is consistent with the PI based secondary 
frequency control \cite[Chapter 9]{Jan_Machowski} 
and \cite[Chapter 4]{Hassan_Bevrani}. 
In PIAC, the weights of the frequency deviation of node $i$ are specified as the inertia $M_i$
and the droop coefficient $D_i$ for the proportional input and 
integral input respectively. The proportional input is used 
to estimate the power stored in the inertia at the transient phase, which is usually 
neglected in the traditional ACE method.
The control gain coefficient $k_r$ can 
be different for each area, which can be tuned 
according to the sensitivity of the control devices in the area. 
By (\ref{eq:multi}), as the synchronous frequency deviation is steered to zero, $P_{ex}$ also converges to the nominal value $P_{ex}^*$. 
Similar as the derivation of (\ref{eq:piac_u}) for PIAC (\ref{eq:PIAC}), we derive for (\ref{eq:mulit}) that 
\begin{eqnarray*}
\hspace{50pt} \dot{u}_r(t)=-k_r\Big(\sum_{i\in\mathcal{V}_{A_r}}{P_i}+P_{ex}^*+u_r(t)\Big)
\end{eqnarray*}
which indicates that \emph{the controllers in area $A_r$ only respond to the power imbalance $\sum_{i\in\mathcal{V}_{A_r}}{P_i}$}. 
Hence, in the network, the control actions of all the areas can be done in \emph{an asynchronous way
where each area can balance the local power supply-demand 
at any time according to the availability of the devices}.

In particular, PIAC becomes a decentralized control method if each node is seen as a single area and $\mathcal{V}_{P}=\emptyset$, 
$\mathcal{V}_K=\mathcal{V}_M\cup\mathcal{V}_F$, i.e., for all $i\in\mathcal{V}_M\cup\mathcal{V}_F$,
\begin{subequations}
 \begin{align*}
 \dot{\eta}_i&=D_i\omega_i+\sum_{j\in\mathcal{V}}B_{ij}\sin{(\theta_{ij})}-\sum_{j\in\mathcal{V}}B_{ij}\sin{(\theta_{ij}^*)},\\
  u_i&=-k_i(M_i\omega+\eta_i),
 \end{align*}
\end{subequations}
where  $\{\theta_{i}^*,~i\in\mathcal{V}\}$ are 
the steady state calculated by tertiary control. 
Since $u_i$ is tracking $-P_i+\sum_{j\in\mathcal{V}}{B_{ij}\sin{(\theta_{ij}^*})}$, each node compensates the power imbalance locally and 
the control actions of the nodes are irrelevant to each other. However, the control cost is not optimized by 
this decentralized control method. 

\section{Stability analysis of PIAC}\label{Sec:stability}

In this section, we analyze the stability of PIAC with the Lyapunov-LaSalle approach as in \cite{DePersis2016,Trip2016240}.
The stability proof makes use of Theorem A.1 stated in the Appendix. 

As indicated in subsection \ref{subsec:multi_description}, the control actions of the areas 
are decoupled in the multi-area control network. So we only need to prove the stability of PIAC implemented in a single-area network.
Extension to multi-area control networks then follows immediately. With the control law (\ref{eq:PIAC}), 
the closed-loop system of PIAC is 
\begin{subequations}\label{PIAC:theta}
 \begin{align}
  \dot{\theta}_i&=\omega_i,~i\in\mathcal{V}_M\cup\mathcal{V}_F,\label{eq:cls1}\\
M_{i}\dot{\omega}_{i}&=P_{i}-D_{i}\omega_{i}-\sum_{j\in \mathcal V}{B_{ij}\sin{(\theta_{ij})}}+{J'}_{i}^{-1}(\lambda),~ 
i\in \mathcal{V}_M ,\\
0&=P_{i}-D_{i}\omega_{i}-\sum_{j\in \mathcal V}{B_{ij}\sin{(\theta_{ij})}}+{J'}_{i}^{-1}(\lambda), ~i\in \mathcal{V}_F ,\\
0&=P_{i}-\sum_{j\in \mathcal V}{B_{ij}\sin{(\theta_{ij})}},~ i\in \mathcal{V}_P , \\
\dot{\eta}&=\sum_{i\in \mathcal{V}_M \cup \mathcal{V}_F }D_i\omega_i,\label{eq:integral}\\
0&=\sum_{i\in\mathcal{V}_K}{{J'}_i^{-1}(\lambda)}+k(\sum_{i\in \mathcal{V}_M }M_i\omega_i+\eta),\\
0&={J'}_{i}^{-1}(\lambda),i~\notin\mathcal{V}_K,\label{eq:clsn} 
 \end{align}
\end{subequations}
where $\theta_{ij}=\theta_i-\theta_j$ is the phase angle differences between 
the nodes connected by a transmission line. 

We denote the angles in the sets $\mathcal{V}_M,\mathcal{V}_F,\mathcal{V}_P$ by column vectors
$\theta_M,\theta_F,\theta_P$, the frequency deviations by column vectors
$\omega_M,\omega_F,\omega_P$,  the angles in $\mathcal{V}$ by $\theta=(\theta_M^T,\theta_F^T,\theta_P^T)^T$,
and the frequency deviations by $\omega=(\omega_M^T,\omega_F^T,\omega_P^T)^T$. 

Note that the closed-loop system may not have a synchronous state if the power injections $\{P_i,i\in\mathcal{V}\}$ 
are much larger than the line capacity $B_{ij}$. For more details on the synchronous state of the power system, we refer to \cite{Xi2016,Dorfler2012}. 
Therefore, we make Assumption \ref{assumption4}.

\begin{assumption}\label{assumption4}
For the closed-loop system (\ref{PIAC:theta}), there exists a synchronous state
$(\theta^*, \omega^*,\eta^*,\lambda^*)\in \mathbb{R}^n \times \mathbb{R}^{n}\times \mathbb{R}\times \mathbb{R}$ with $\omega^*_i=\omega_{syn}$   
and $$\theta^*(t)\in\Theta=\{\theta\in \mathbb{R}^n\lvert|\theta_i-\theta_j|<\frac{\pi}{2}, ~\forall (i,j)\in\mathcal{E}\}.$$ 
\end{assumption} 
The condition $\theta^*\in \Theta$ is commonly referred to as a security constraint \cite{DePersis2016} in power 
system analysis. It can be satisfied by reserving some margin of power flow when calculating the operating point in 
tertiary control \cite{Ilic2000}.

Since the power flows $\{B_{ij}\sin(\theta_i-\theta_j), ~(i,j)\in\mathcal{E}\}$ only depend on the angle differences and 
the angles can be expressed relative to a reference node, we choose a reference angle, i.e., $\theta_1$, in $\mathcal{V}_M$
and introduce the new variables
\begin{eqnarray*}
 \hspace{50pt} \varphi_i=\theta_i-\theta_1, ~i=1,2,\cdots, n
\end{eqnarray*}
which yields $\dot{\varphi_i}=\omega_i-\omega_1$. Note that 
$\varphi_1=0, \dot{\varphi}_1=0$ for all $t>0$. With 
$\omega_i=\dot{\varphi_i}+\omega_1$ for $i\in\mathcal{V}_F$, the closed loop system (\ref{PIAC:theta})
can be written in the DAE form as (\ref{eq:DAEtt}) in the Appendix,
\begin{subequations}\label{PIAC:varphi}
 \begin{align}
  \dot{\varphi}_i&=\omega_i-\omega_1,~i\in\mathcal{V}_M\cup\mathcal{V}_F,\label{eq:nc1}\\
\hspace{-10pt}
M_{i}\dot{\omega}_{i}&=P_{i}-D_{i}\omega_{i}-\sum_{j\in \mathcal V}{B_{ij}\sin{(\varphi_{ij})}}+{J'}_{i}^{-1}(\lambda), 
i\in \mathcal{V}_M ,\label{eq:nc2}\\
D_i\dot{\varphi}_i&=P_{i}-D_{i}\omega_{1}-\sum_{j\in \mathcal V}{B_{ij}\sin{(\varphi_{ij})}}+{J'}_{i}^{-1}
(\lambda), i\in \mathcal{V}_F ,\label{eq:nc3}\\
0&=P_{i}-\sum_{j\in \mathcal V}{B_{ij}\sin{(\varphi_{ij})}}, i\in \mathcal{V}_P,\label{eq:nc4} \\
\dot{\eta}&=\sum_{i\in \mathcal{V}_M \cup \mathcal{V}_F }D_i\omega_i,\label{eq:nc5}\\
0&=\sum_{i\in\mathcal{V}_K}{{J'}_i^{-1}(\lambda)}+k(\sum_{i\in \mathcal{V}_M }M_i\omega_i+\eta),\label{eq:nc6}
 \end{align}
\end{subequations}
where $\varphi_{ij}=\varphi_i-\varphi_j$, and the equations (\ref{eq:nc1}-\ref{eq:nc4}) are from the power system and (\ref{eq:nc5}-\ref{eq:nc6}) 
from the controllers. We next recast the system (\ref{PIAC:varphi}) into the form of the DAE system (\ref{eq:DAEtt}), the state variables are 
$x=(\varphi_{M},\varphi_{F},\omega_M,\eta)\in \mathbb{R}^{n_M-1}\times \mathbb{R}^{n_F}\times \mathbb{R}^{n_M}\times \mathbb{R}$, 
the algebraic variables are $y=(\varphi_{P},\lambda)\in \mathbb{R}^{n_P}\times \mathbb{R}$,
the differential equations are (\ref{eq:nc1}-\ref{eq:nc3},\ref{eq:nc5}) and the algebraic equations are (\ref{eq:nc4}, \ref{eq:nc6}). 
 Here $\varphi_{M}$ is with the components $\{\varphi_i,i\in \mathcal{V}_M\}$ besides $\varphi_1$ which is a constant, $\varphi_{F}$ with components 
$\{\varphi_i,i\in \mathcal{V}_F\}$, and $\varphi_{P}$ with components 
$\{\varphi_i,i\in \mathcal{V}_P\}$. 
Note that the variables $\{\omega_i, i\in\mathcal{V}_F\}$ are not included 
into the state variable or algebraic variables since the terms $\{D_i\omega_i,i\in\mathcal{V}_F\}$ in (\ref{eq:nc5}) can be replaced by 
$\{P_i-\sum_{j\in \mathcal V}{B_{ij}\sin{(\varphi_{ij})}}+{J'}_{i}^{-1}(\lambda),i\in\mathcal{V}_F\}$
. (\ref{eq:clsn}) is neglected since it is irrelevant to the following stability analysis.

When mapping $\theta$ to the coordinate of $\varphi$, Assumption \ref{assumption4}
yields $$\varphi\in\Phi=\{\varphi\in \mathbb{R}^n||\varphi_i-\varphi_j|<\frac{\pi}{2},~\forall (i,j)\in\mathcal{E}, ~\text{and}~ \varphi_1=0\}.$$
We remark
that each equilibrium state of (\ref{PIAC:varphi}) corresponds to a synchronous state of (\ref{PIAC:theta}). 
In the new coordinate, 
we have the following theorem for the stability of the system (\ref{PIAC:varphi}). 
\begin{theorem}\label{thm2}
 If the assumptions \ref{assumption1}-\ref{assumption4} hold, for the system (\ref{PIAC:varphi}), then
 \begin{enumerate}[(a)]
 \item there exists an unique equilibrium state
 $z^*=(\varphi^*,\omega_M^*,\eta^*,\lambda^*)\in\Psi$ where $\Psi=\Phi\times\mathbb{R}^{n_M}\times\mathbb{R}\times \mathbb{R}$. 
 \item there exists a domain $\Psi^d\subset\Psi$ such that 
 for any initial state $z^0\in \Psi^d$ that satisfies
 the algebraic equations (\ref{eq:nc4}) and (\ref{eq:nc6}), the state trajectory converges to the unique equilibrium state $z^*\in\Psi$.
 \end{enumerate}
\end{theorem}
Note that the cost functions are not required to be scaled for the local asymptotically stable of PIAC as assumed in \cite{Dorfler2016}, and 
the size of the attraction domain of the equilibrium state $z^*$ is not determined in Theorem \ref{thm2}
which states the stability of the PIAC method.
The proof of Theorem \ref{thm2} is based on the Lyapunov/LaSalle stability criterion as in Theorem A.1. 
In the following, we first present the verification of the Assumption A.1 and A.2 in the Appendix for the DAE system (\ref{PIAC:varphi}), 
then prove the stability of (\ref{PIAC:varphi}) by designing a Lyapunov function as $V(x,y)$ in Theorem A.1. 
Lemma \ref{lem1} states that (\ref{PIAC:varphi}) possesses an 
equilibrium state, which verifies Assumption A.1,
and Lemma \ref{lem2} claims the regularity of the algebraic equation (\ref{eq:nc4}, \ref{eq:nc6}), which verifies Assumption A.2.

\begin{lemma}\label{lem1}
There exists at most one equilibrium state  $z^*=(\varphi^*, \omega^*_M,\eta^*,\lambda^*)$ of 
the system (\ref{PIAC:varphi}) such that 
$z^*\in\Psi$ and
\begin{subequations}
 \begin{align}
  \omega^*_i&=0, ~i\in\mathcal{V}_M,\\
 \sum_{i\in\mathcal{V}}{P_i}+\sum_{i\in\mathcal{V}_K}{{J'}^{-1}_i(\lambda^*)}&=0.\label{eq:equilibrium2}\\
 \sum_{i\in\mathcal{V}}{P_i}-k\eta^*&=0. \label{eq:equilibrium3}
 \end{align}
\end{subequations}

\end{lemma}
\emph{Proof:} At the synchronous state, from (\ref{eq:piac_u}) and (\ref{eq:PIACb}) it follows that 
\begin{eqnarray}
 ~~~~~~~~~~~~~~~k(\sum_{i\in \mathcal{V}_M }M_i\omega^*_i+\eta^*)=\sum_{i\in\mathcal{V}}{P_i}.\label{eq:eql1}
\end{eqnarray}
Substitution into the algebraic equation (\ref{eq:nc6}) yields (\ref{eq:equilibrium2}). 
Substitution of (\ref{eq:equilibrium2}) into (\ref{eq:sync}) with $u_i={J'}^{-1}_i(\lambda^*)$, we obtain $\omega_{syn}=0$ which yields 
$\{\omega^*_i=0,i\in\mathcal{V}_M\}$. Hence this and (\ref{eq:eql1}) yield
$$\sum_{i\in\mathcal{V}_K}{{J'}^{-1}_i(\lambda^*)}+k\eta^*=0.$$
which leads to (\ref{eq:equilibrium3}). 
It follows from \cite{Araposthatis1981,skar_uniqueness_equilibrium} that the system (\ref{PIAC:theta}) has at most one power flow solution such that $\theta\in\Theta$. Hence 
there exists at most one equilibrium for the system (\ref{PIAC:varphi}) that satisfies $\varphi\in\Phi$. 
\hfill $\Box$

With respect to the regularity of the algebraic equations (\ref{eq:nc4}, \ref{eq:nc6}), we derive the following lemma.
\begin{lemma}\label{lem2}
For any $\varphi\in\Phi$ and strictly 
convex functions of $\{J_i(u_i),i\in\mathcal{V}_K\}$ in the optimization problem (\ref{eq:optimal1}), 
the algebraic equations (\ref{eq:nc4}, \ref{eq:nc6}) 
are regular.
\end{lemma}
\emph{Proof:} Since (\ref{eq:nc4}) and (\ref{eq:nc6}) are independent algebraic equations 
with respect to $\varphi_P$ and $\lambda$ respectively, the regularity of each of them 
is proven separately.  

First, we prove the regularity of (\ref{eq:nc4}) by showing that its Jacobian is a principle minor of the 
Laplacian matrix of a weighted network. In the coordination of $\theta$, we define 
function $$\overline{U}(\theta)=\sum_{(i,j)\in\mathcal{E}}{B_{ij}(1-\cos{(\theta_i-\theta_j))}}.$$ The Hessian matrix 
of $\overline{U}(\theta)$ is 
\begin{eqnarray*}
 \overline{L}=\left(\begin{array}{cccc} \overline{B}_{11} & -B_{12}\cos{(\theta_{12})} & \dots& -B_{1n}\cos{(\theta_{1n})}  \\
                            -B_{21}\cos{(\theta_{12})} & \overline{B}_{22} & \dots& -B_{2n}\cos{(\theta_{2n})}  \\
                            \vdots &\vdots&\ddots&\vdots \\
                            -B_{n1}\cos{(\theta_{n1})} & -B_{n2}\cos{(\theta_{n2})} & \dots & \overline{B}_{nn}
                       \end{array}\right), 
\end{eqnarray*}

where $\overline{B}_{ii}=\sum_{j\in\mathcal{V}}{B_{ij}\cos{(\theta_{ij})}}$ and $\theta_{ij}=\theta_i-\theta_j$. $\overline{L}$ is the Laplacian of 
the undirected graph $G$ defined in section {\ref{Sec:model}} with positive line weights 
$B_{ij}\cos{(\theta_i-\theta_j)}$. Hence $\overline{L}$ is semi-positive definite and all its principle minors
are nonsingular\cite[Theorem 9.6.1]{Brualdi}. In the coordination of $\varphi$, we define function 
\begin{eqnarray}\label{eq:def_U}
     U(\varphi)=\sum_{(i,j)\in\mathcal{E}}{B_{ij}(1-\cos{(\varphi_i-\varphi_j))}},  
       \end{eqnarray}
with Hessian matrix 
\begin{small}
\begin{eqnarray}\label{eq:Hess}
\hspace{-4pt}
 L=\left(\begin{array}{cccc} B_{22} & -B_{23}\cos{(\varphi_{23})} & \dots& -B_{2n}\cos{(\varphi_{2n})}  \\
                            -B_{32}\cos{(\varphi_{32})} & B_{33} & \dots& -B_{3n}\cos{(\varphi_{3n})}  \\
                            \vdots &\vdots&\ddots&\vdots \\
                            -B_{n2}\cos{(\varphi_{n2})} & -B_{n2}\cos{(\varphi_{n3})} & \dots & B_{nn}
                       \end{array}\right) 
\end{eqnarray}
\end{small}
where $B_{ii}=\sum_{j\in\mathcal{V}}{B_{ij}\cos{(\varphi_{ij})}}$, $\varphi_{ij}=\varphi_i-\varphi_j$, and $\varphi_1=0$. Note that 
$B_{ij}\cos{(\varphi_{ij})}=B_{ij}\cos{(\theta_{ij})}$, thus $L$ is a principle minor of $\overline{L}$ and is nonsingular. Hence  
the Jacobian of (\ref{eq:nc4}) with respect to $\varphi_{P}$, which is a principle minor of $L$, is nonsingular. 

Second, because $J_i$ is strictly convex by Assumption \ref{assumption2} such that $J''_i>0$, we obtain $({J'}_i^{-1})'=\frac{1}{J''_i}>0$ which yields
$\big(\sum_{i\in\mathcal{V}_K}{{J'}_i^{-1}(\lambda)}\big)'>0$. Hence (\ref{eq:nc6}) is nonsingular. \hfill $\Box$

What follows is the proof of Theorem \ref{thm2} with the Lyapunov/LaSalle stability criterion.   

\emph{Proof of Theorem \ref{thm2}:} 
Lemma \ref{lem1} and Assumption \ref{assumption4} states that the equilibrium state $z^*=(\varphi^*,\omega_M^*,\eta^*,\lambda^*)\in\Psi$ is unique.  
The proof of the statement (b) in Theorem \ref{thm2} is based on Theorem A.1.  Consider an incremental Lyapunov function candidate, 
\begin{eqnarray}
V(\varphi,\omega_M,\eta,\lambda)=V_1+\alpha V_2+V_3, \label{Lyapunov_function}
\end{eqnarray}
where $V_1$ is the classical energy-based function \cite{DePersis2016}, 
\begin{eqnarray*}
   && V_1(\varphi, \omega_M)=U(\varphi)-U(\varphi^*)-\nabla_{\varphi}U(\varphi^*)(\varphi-\varphi^*)\\
  &&~~~~~~~~~~~~~~~~+\frac{1}{2}\omega_M^TM_M\omega_M,
\end{eqnarray*}

and $V_2, V_3$ are positive definite functions
\begin{eqnarray*}
&&~~V_2(\lambda)\quad=\frac{1}{2}(\sum_{i\in\mathcal{V}_K}{{J'}_i^{-1}(\lambda)}-\sum_{i\in\mathcal{V}_K}{{J'}_i^{-1}(\lambda^*)})^2,\\
&&V_3(\omega_M,\eta)=\frac{k^2}{2}(\sum_{i\in \mathcal{V}_M }M_i\omega_i+\eta-\eta^*)^2.
\end{eqnarray*}

Note that the definition of $U(\varphi)$ is in (\ref{eq:def_U}) and $V_2=V_3$ by (\ref{eq:nc6}). $V_3$ is introduced to involve 
 the state variable $\eta$ into the Lyapunov function.  

First, we prove that $\dot{V}\leq 0$. From the dynamic system (\ref{PIAC:varphi}) and the definition of $V_1$, 
it yields that 
\begin{small}
\begin{eqnarray}
 \dot{V}_1=-\sum_{i\in\mathcal{V}_M\cup\mathcal{V}_F}{D_i\omega_i^2}+
 \sum_{i\in\mathcal{V}_K}\omega_i({J'}_i^{-1}(\lambda)-{J'}_i^{-1}(\lambda^*)),\label{Lyp:1}
\end{eqnarray}
\end{small}
 by ({\ref{eq:nc6}}), we derive
\begin{small}
\begin{subequations}
 \begin{align}
 \hspace{-15pt}
 \nonumber
 &\dot{V}_2=-k\Big[\sum_{i\in\mathcal{V}_K}{{J'}_i^{-1}(\lambda)}-\sum_{i\in\mathcal{V}_K}{{J'}_i^{-1}(\lambda^*)}\Big]
 \Big[\sum_{i\in \mathcal{V}_M }M_i\dot{\omega}_i
 +\dot{\eta}\Big]\\
 \nonumber
  \hspace{-15pt}
&~~\text{by summing (\ref{eq:nc2}-\ref{eq:nc5})}\\
 \nonumber
  \hspace{-15pt}
 =&-k\big[\sum_{i\in\mathcal{V}_K}{{J'}_i^{-1}(\lambda)}-\sum_{i\in\mathcal{V}_K}{{J'}_i^{-1}(\lambda^*)}\big]\big[\sum_{i\in\mathcal{V}}{P_i}
 +\sum_{i\in\mathcal{V}_K}{{J'}_i^{-1}(\lambda)}\big]\\
  \nonumber
   \hspace{-15pt}
 &~~\text{by (\ref{eq:equilibrium2})}\\
 =&-k\Big[\sum_{i\in\mathcal{V}_K}{{J'}_i^{-1}(\lambda)}-\sum_{i\in\mathcal{V}_K}{{J'}_i^{-1}(\lambda^*)}\Big]^2,\label{eq:dotV3}\\
  \nonumber
   \hspace{-15pt}
 &~~\text{by expanding the quadratic}\\
 \nonumber
  \hspace{-15pt}
 =&-k\sum_{i\in\mathcal{V}_K}\big[{J'}_i^{-1}(\lambda)-{J'}_i^{-1}(\lambda^*)\big]^2\\
 \nonumber
  \hspace{-15pt}
 &-2k\sum_{i\neq j,i,j\in\mathcal{V}_K}
 {\Big[{J'}_i^{-1}(\lambda)-{J'}_i^{-1}(\lambda^*)\Big]\Big[{J'}_j^{-1}(\lambda)-{J'}_j^{-1}(\lambda^*)\Big]}\\
  \nonumber
   \hspace{-15pt}
 &~~\text{by Newton-Leibniz formula}\\
 \nonumber
  \hspace{-15pt}
 =& -k\sum_{i\in\mathcal{V}_K}\big[{J'}_i^{-1}(\lambda)-{J'}_i^{-1}(\lambda^*)\big]^2\\
   \hspace{-15pt}
 &-2k\sum_{i\neq j,i,j\in\mathcal{V}_K}{\Big[\int_{\lambda^*}^{\lambda}({J'}_i^{-1})'ds\Big]
 \Big[\int_{\lambda^*}^{\lambda}({J'}_j^{-1})'ds \Big]}.\label{Lyp:2}
\end{align}
\end{subequations}
\end{small}
and by $V_3=V_2$ and (\ref{eq:dotV3}), we obtain 
\begin{eqnarray}
 \dot{V}_3=-k\Big[\sum_{i\in\mathcal{V}_K}({J'}_i^{-1}(\lambda)-{J'}_i^{-1}(\lambda^*))\Big]^2. \label{Lyp:3}
\end{eqnarray}
Hence,  with (\ref{Lyp:1},\ref{Lyp:2},\ref{Lyp:3}), we derive
\begin{small}
\begin{subequations}
 \begin{align*}
  &\dot{V}=-\sum_{i\in\mathcal{V}_M\cup\mathcal{V}_F}{D_i\omega_i^2}+
 \sum_{i\in\mathcal{V}_K}\omega_i\big[{J'}_i^{-1}(\lambda)-{J'}_i^{-1}(\lambda^*)\big]\\
 &\quad
 -k\alpha\sum_{i\in\mathcal{V}_K}\big[{J'}_i^{-1}(\lambda)-{J'}_i^{-1}(\lambda^*)\big]^2  \\
 &\quad
 -k\alpha\sum_{i\neq j,i,j\in\mathcal{V}_K}{\Big[\int_{\lambda^*}^{\lambda}({J'}_i^{-1})'ds\Big]
 \Big[\int_{\lambda^*}^{\lambda}({J'}_j^{-1})'ds \Big]}\\
 &\quad
 -k\Big[\sum_{i\in\mathcal{V}_K}{J'}_i^{-1}(\lambda)-\sum_{i\in\mathcal{V}_K}{J'}_i^{-1}(\lambda^*)\Big]^2\\
 &=-\sum_{i\in\mathcal{V}_K}{\frac{D_i}{2}\omega_i^2}-
 \sum_{i\in\mathcal{V}_K}{\frac{D_i}{2}\Big[\omega_i-\frac{{J'}_i^{-1}(\lambda)-{J'}_i^{-1}(\lambda^*)}{D_i}\Big]^2}\\
 &\quad-\sum_{i\in\mathcal{V}_K}\big[k\alpha-\frac{1}{2D_i}\big]\big[{J'}_i^{-1}(\lambda)-{J'}_i^{-1}(\lambda^*)\big]^2
  \\
 &\quad
 -2k\alpha\sum_{i\neq j,i,j\in\mathcal{V}_K}{\Big[\int_{\lambda^*}^{\lambda}({J'}_i^{-1})'ds\Big]
 \Big[\int_{\lambda^*}^{\lambda}({J'}_j^{-1})'ds \Big]}\\
 &\quad
 -k\Big[\sum_{i\in\mathcal{V}_K}{J'}_i^{-1}(\lambda)-\sum_{i\in\mathcal{V}_K}{J'}_i^{-1}(\lambda^*)\Big]^2-
 \sum_{i\in\mathcal{V}_M\cup\mathcal{V}_F/\mathcal{V}_K}{D_i\omega_i^2}
 \end{align*}
\end{subequations}
\end{small}
where the equation 
$$\sum_{i\in\mathcal{V}_M\cup\mathcal{V}_F}{D_i\omega_i^2}=\sum_{i\in\mathcal{V}_K}{D_i\omega_i^2}+
\sum_{i\in\mathcal{V}_M\cup\mathcal{V}_F\backslash\mathcal{V}_K}{D_i\omega_i^2}$$ is used 
due to the fact that $\mathcal{V}_K\subseteq\mathcal{V}_M\cup\mathcal{V}_F$. 

Since $J_i(u_i)$ is strictly convex and $({J'}_i^{-1})'=\frac{1}{{J''}_{i}}>0$, we derive 
$$k\alpha\sum_{i\neq j,i,j\in\mathcal{V}_K}{\Big[\int_{\lambda^*}^{\lambda}({J'}_i^{-1})'ds\Big]
 \Big[\int_{\lambda^*}^{\lambda}({J'}_j^{-1})'ds \Big]}>0,$$
 Thus by setting $\alpha>\frac{1}{kD_i}$, we obtain $\dot{V}\leq 0$. 

Second, we prove that $z^*=(\varphi^*,\omega^*_M,\eta^*,\lambda^*)$ is a 
strict minimum of $V(\varphi,\omega_M,\eta,\lambda)$ such 
that $\nabla V|_{z^*}=0$ and $\nabla^2 V|_{z^*}>0$. 
It can be easily verified that $V|_{z^*}=0$ and 
\begin{eqnarray*}
 \nabla V|_{z^*}=\text{col}(\nabla_\varphi V, \nabla_{\omega_M} V, \nabla_\eta V,\nabla_\lambda V)|_{z^*}=0\in \mathbb{R}^{n+n_M+1}
\end{eqnarray*}
where
\begin{eqnarray*}
\hspace{-5pt}
&&\nabla_\varphi V=\nabla_\varphi U-\nabla_{\varphi} U^*,\\
\hspace{-5pt}
&&\nabla_{\omega_M} V=M\omega_M+M \big(k^2(\sum_{i\in \mathcal{V}_M }M_i\omega_i+\eta-\eta^*))1_{n_M}\big),\\
\hspace{-10pt}
&&\nabla_\eta V=k^2(\sum_{i\in \mathcal{V}_M }M_i\omega_i+\eta-\eta^*)),\\
\hspace{-5pt}
&&\nabla_\lambda V=\alpha\Big(\sum_{i\in\mathcal{V}_K}{{J'}_i^{-1}(\lambda)}-\sum_{i\in\mathcal{V}_K}{{J'}_i^{-1}(\lambda^*)}\Big)
\Big(\sum_{i\in\mathcal{V}_K}{{J'}_i^{-1}(\lambda)}\Big)'.
\end{eqnarray*}
Here $k^2(\sum_{i\in \mathcal{V}_M }M_i\omega_i+\eta-\eta^*)1_{n_M}$ is a vector with all components equal to
$k^2(\sum_{i\in \mathcal{V}_M }M_i\omega_i+\eta-\eta^*))$.
The Hessian matrix of $V$ is
\begin{eqnarray*}
\hspace{60pt} \nabla^2 V|_{z^*}=\text{blkdiag}(L,H,\Lambda),
\end{eqnarray*}
which is a block diagonal matrix with block matrices $L$,$H$, and $\Lambda$. 
$L$ is positive definite by (\ref{eq:Hess}),
$$\Lambda=\alpha \Big(\big(\sum_{i\in\mathcal{V}_K}{{J'}_i^{-1}(\lambda^*)}\big)'\Big)^2>0$$
which is a scalar value, and $H$ is the Hessian matrix of the function 
\begin{eqnarray*}
 \hspace{20pt}\overline{V}=\frac{1}{2}\omega_M^TM_M\omega_M+\frac{k^2}{2}\Big(\sum_{i\in \mathcal{V}_M }M_i\omega_i+\eta-\eta^*\Big)^2
\end{eqnarray*}
which is positive definite for any $(\omega_M,\eta-\eta^*)$, thus $H$ is positive definite. Hence, we have proven 
that $z^*$ is a strict minimum of $V$. 
 
 Finally, we prove that the invariant set 
$$\Big\{(\varphi,\omega_M,\eta,\lambda)|\dot{V}((\varphi,\omega_M,\eta,\lambda)=0\Big\}$$ contains only the
equilibrium point. $\dot{V}=0$ implies that $\{\omega_i=0,i\in\mathcal{V}_M\cup\mathcal{V}_F\}$. Hence
$\{\varphi_i,i\in\mathcal{V}\}$ are constants. By lemma \ref{lem1}, there is at most one equilibrium with 
$\varphi\in\Phi$. In this case, $z^*$ is the only one equilibrium in the neighborhood of $z^*$ , i.e.,
$\Psi^d=\{(\varphi,\omega_M,\eta,\lambda)|V(\varphi,\omega_M,\eta,\lambda)\leq c, \varphi\in\Phi\}$ for some 
$c>0$. Hence with any initial state $z^0$ that satisfies the algebraic equations (\ref{eq:nc4}) and (\ref{eq:nc6}), the 
trajectory converges to the equilibrium state $z^*$. \hfill $\Box$

For the multi-area implementation of PIAC, we choose a Lyapunov candidate function as 
\begin{eqnarray*}
 V(\varphi,\omega_M,\eta,\lambda)=V_1+\sum_{A_r}(\alpha V_{2r}+V_{3r}),
\end{eqnarray*}
where  $\lambda=\text{col}(\lambda_r)$ is a column vector consisting of the components $\lambda_r$, 
$V_1$ is defined as in (\ref{Lyapunov_function}) and $V_{2r}$ and $V_{3r}$ are defined for area $A_r$ as 
\begin{subequations}
 \begin{align*}
  V_{2r}&=\frac{1}{2}\big(\sum_{i\in\mathcal{V}_{K_r}}{J'}_i^{-1}(\lambda_r)-\sum_{i\in\mathcal{V}_{K_r}}{J'}_i^{-1}(\lambda^*_r)\big),\\
  V_{3r}&=\frac{k^2}{2}\big(\sum_{i\in\mathcal{V}_{M_r}}{M_i\omega_i}+\eta_r-\eta_r^*\big)^2.
 \end{align*}
\end{subequations}
Following the proof of Theorem \ref{thm2}, we can obtain the locally asymptotic stability of PIAC implemented in multi-area control.

\begin{remark}\label{remark2}
The assumptions \ref{assumption1}-\ref{assumption4} are realistic at the same time. 
Assumption \ref{assumption1} and \ref{assumption3} are necessary
for the implementation of PIAC to solve Problem \ref{problem1} and \ref{problem2}. Assumption \ref{assumption2} and \ref{assumption4} are 
general sufficient conditions for the stability of the nonlinear systems (\ref{eq:model1}) controlled by PIAC. 
Assumption \ref{assumption1} and \ref{assumption4} can be guaranteed by tertiary control and Assumption \ref{assumption3} by an
effective communication infrastructure. Assumption \ref{assumption2} usually holds for frequently used convex cost functions, e.g., quadratic 
cost function, where the requirement of scale cost functions in \cite[Assumption 1]{Dorfler2016}
is not needed. 
\end{remark}

 \section{Simulations of the closed-loop system}\label{Sec:experiment}

\begin{figure}[ht]
\centering
 \includegraphics[scale=0.32]{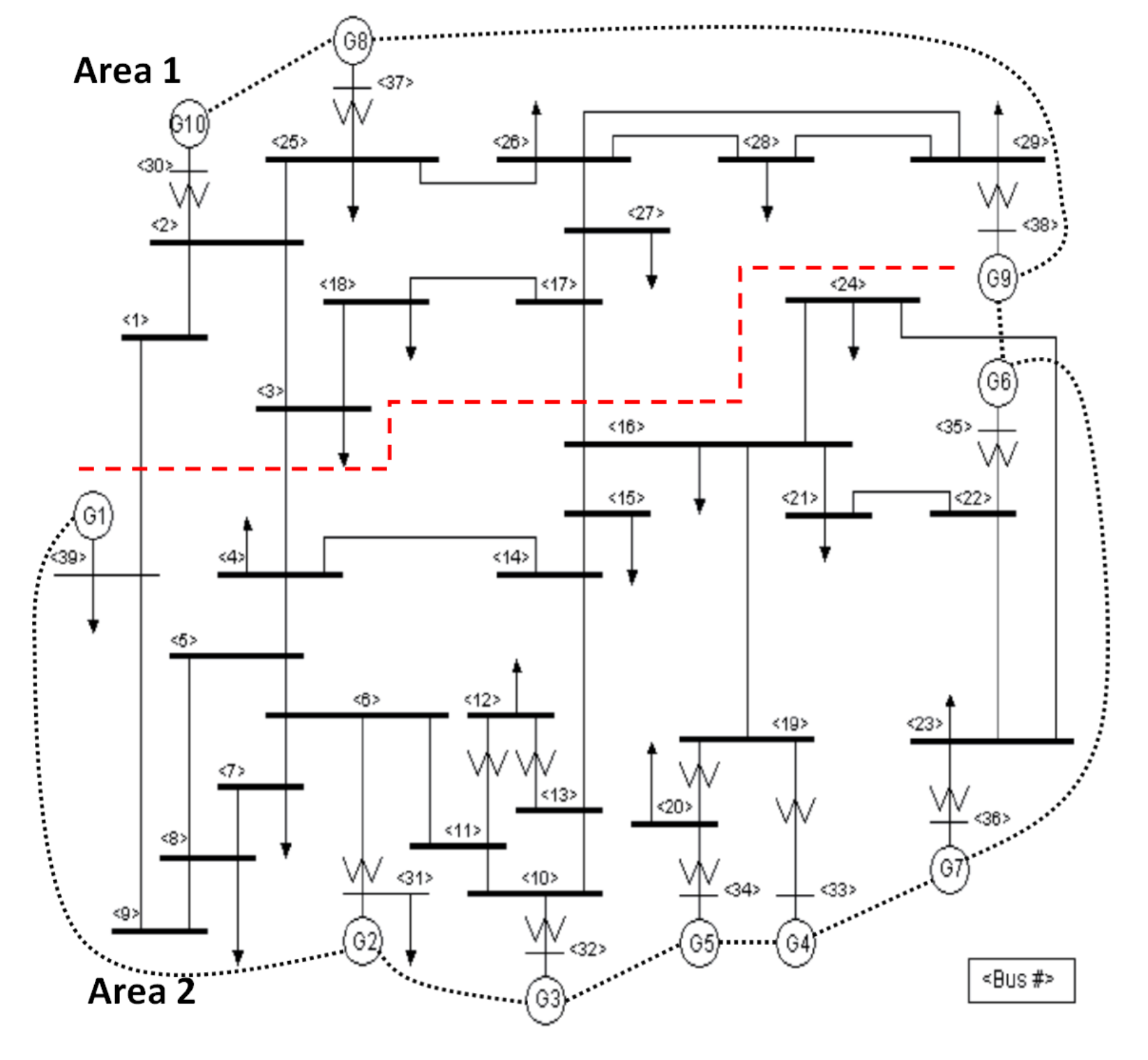}
 \caption{IEEE New England test power system.}
 \label{fig.IEEE39_graph}
\end{figure}

In this section, we evaluate the performance of the PIAC method and compare it with 
with those of the GB, DAI and DecI control laws on the
IEEE New England power grid shown in Fig.~\ref{fig.IEEE39_graph}. The data 
are obtained from \cite{Athay1979}. 
In the test system, there are 10 generators, and 39 buses and it serves a total load of about 6 GW. The voltage at each bus
is a constant which is obtained by power flow calculation with the \emph{Power System Analysis Toolbox} (PSAT) \cite{Milano2008}. 
In the network, 
there actually are 49 nodes, i.e., 10 nodes for the generators, 39 nodes for the buses. Each synchronous machine is connected 
to a bus and its phase angle is rigidly tied to the rotor angle of the bus if the 
voltages of the system are constants, e.g., the classical model of synchronous machines \cite{Ilic2000}. 
We simplify the test system to a system of 39 nodes by considering the generator and the bus as one node. This is 
reasonable because the angles of the synchronous machine and the bus have the same dynamics. 
The 10 generators are in the set $\mathcal{V}_M =\{30,31,32,33,34,35,36,37,38,39\}$ and the other buses are in the 
set $\mathcal{V}_F $ which are assumed to be 
frequency dependent loads. The nodes participating in secondary frequency control are the 10 generators, thus $\mathcal{V}_K =\mathcal{V}_M $. 
The inertias of the generators as stated in \cite{Athay1979} are all divided by 100 in order to obtain 
the desired frequency response as in \cite{Dorfler2016,Zhao2015}. 
The buses in $\mathcal{V}_M \cup \mathcal{V}_F $ and controllers in $\mathcal{V}_K $ are connected by a communication network. 

As in \cite{Dorfler2016,Zhao2015}, we set the droop control coefficient $D_i=1$~(p.u. power/p.u. frequency deviation) for $i\in \mathcal{V}_F  \cup \mathcal{V}_M $
under the power base of 100 MVA and frequency base of 60 Hz with the nominal frequency $f^*=60$ Hz, and 
we choose the quadratic cost function  $J_i(u_i)=\frac{u_i^2}{2a_i},~i\in\mathcal{V}_K$.
The economic dispatch coefficients $a_i$ are generated randomly with a uniform distribution on $(0,1)$. It
can be easily verified that the quadratic cost functions are all strictly convex.

In the simulations, the system is initially at a supply-demand balanced state with the nominal frequency. 
At time $t=0.5$ second, a step-wise increase of 33 MW of the loads at each of the buses 4, 12, and 20, amounting to a total power
imbalance of 99 MW, causes the frequency to 
drop below the nominal frequency. The loads at other nodes do not change.  

We conduct the simulations with the open source software PSAT and use the Euler-Forward 
method to discretize the ordinary differential equations and the Newton-Raphson method to solve the nonlinear system.
We first evaluate the performance of PIAC on the network which is assumed as a single area.
We compare the PIAC method with
the GB, DAI and DecI control laws. 
For illustrations of the overshoot phenomena of other control laws, we refer to the simulation results of the 
literature published recently, e.g., \cite{EAGC,Trip2016240}.
Then we implement the PIAC method in the network separated into two areas and show that the control actions of 
the areas are totally decoupled by PIAC as described in the subsection \ref{subsec:multi_description}.

\subsection{Numerical results of the Single-area implementation}\label{subsec:numerical_single}

\begin{table}[ht]
\centering
\small Table 1. The control parameters\\
 \begin{tabular}{ |c|c|c|c|c|c| }
  \hline
  \multicolumn{1}{|c|}{PIAC}&\multicolumn{2}{|c|}{GB}&\multicolumn{2}{|c|}{DAI}&\multicolumn{1}{|c|}{DecI}\\
  \hline
  $k$ & $k_{GB}$&$C_i$ &$k_i$&$w_{ij}$&$k_i$\\
  \hline
  5 & 60 &$\frac{1}{39}$&50&-20&50\\
  \hline
\end{tabular}
\end{table}

\begin{figure}[ht]
 \begin{center}
  \includegraphics[scale=1.15]{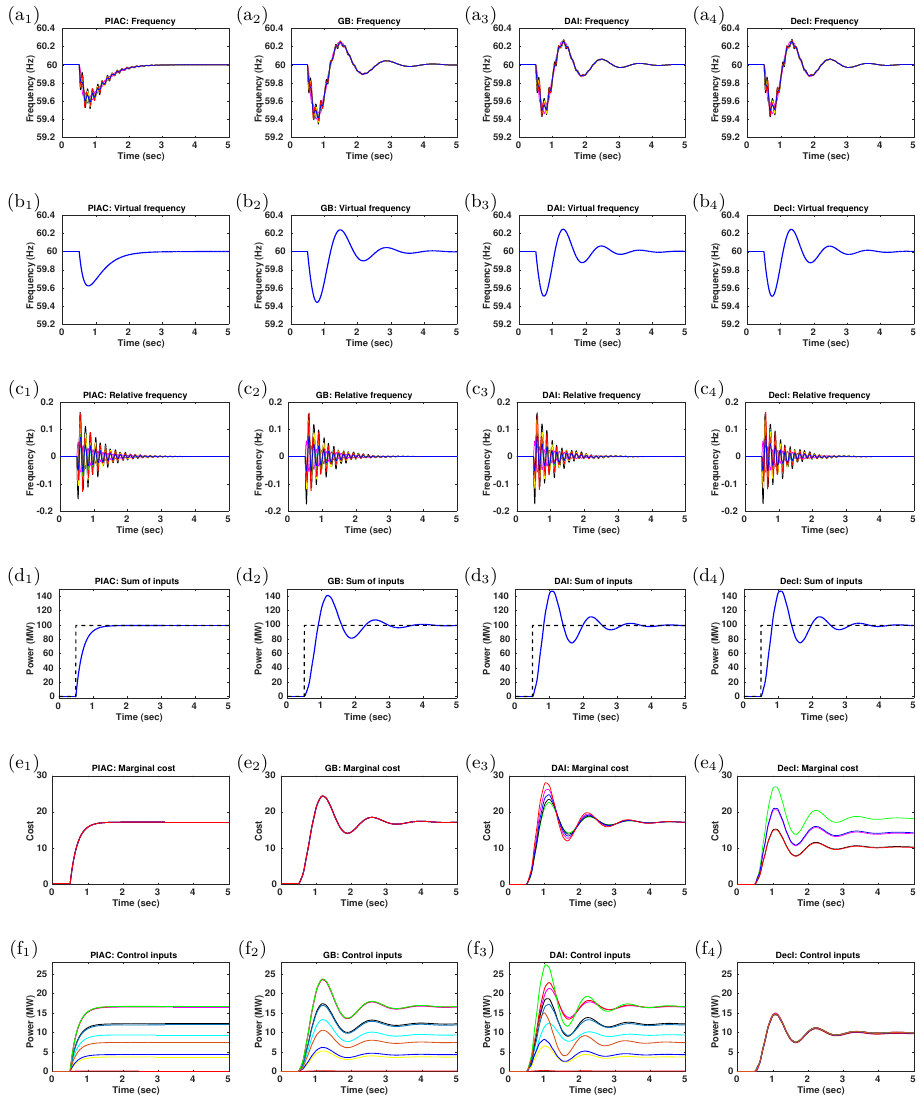}
  \caption{{Comparison of the PIAC, GB, DAI, DecI control laws (single area implementation).
  In ($\text{d}_1$-$\text{d}_4$), the black dashed lines denote the power imbalance of the network.
  Only the marginal costs of the generators $\{2,3,4,5,7\}$ are shown in (e$_3$-e$_4$).}}
  \label{fig.frq}
  \end{center}
  \end{figure}

In this subsection, the network is seen as a single area. 
The parameters for PIAC, GB, DAI and DecI are listed in Table 1. For the DAI method, 
the communication network is a weighted 
network as shown in Fig.~\ref{fig.IEEE39_graph} by the black dotted lines. 
The weight $w_{ij}$ for line $(i,j)$ is set as in Table 1 and $w_{ii}=-\sum_{(i,j)}w_{ij}$. 
The settings of these 
gain coefficients are chosen for a fair comparison in such a way that  
the slopes of the total control inputs, which reflect the required response time of the actuators, 
show close similarity (see Fig.~\ref{fig.frq}d$_1$-\ref{fig.frq}d$_4$).

Fig.~\ref{fig.frq} shows the comparison 
of the performances between the four control laws. Fig.~\ref{fig.frq}a-\ref{fig.frq}d show 
the responses of frequencies $\overline{\omega}_i=\omega_i+f^*$ for all $i\in\mathcal{V}_M$, virtual \emph{abstract
frequency} $\overline{\omega}_s=\omega_s+f^*$, relative frequency $\{\omega_i-\omega_s,~i\in\mathcal{V}_M\}$
and control input $u_s$ respectively. The latter three illustrate 
the three-subprocesses decomposed from the dynamics of the system (1). 
Here, the responses 
of $\omega_s$ are obtained from (10) with $u_s$ as the total amount of 
control inputs of the PIAC and GB method respectively. 
It can be observed from Fig.~\ref{fig.frq}a$_1$ and \ref{fig.frq}a$_4$
that all the control laws can restore the nominal frequency. However,
the frequency deviation under the PIAC method is much smaller than under the other three control laws
which introduce 
extra oscillations to the frequency. This is because the sum of control inputs of 
the GB, DAI, DecI method overshoot the desired value as Fig.~\ref{fig.frq}d$_2$-\ref{fig.frq}d$_4$
show, while the one of the PIAC method 
converges exponentially as Fig.~\ref{fig.frq}$\text{d}_1$ shows. Because of the overshoot, the GB, DAI and DecI control laws
require 
a maximum mechanical power input of about 140 MW from the 10 generators after the disturbance while the PIAC method only 
requires 99 MW in this simulation.
This scenario is also well reflected in the response of 
the virtual frequency $\overline{\omega}_s$ as shown in Fig.\ref{fig.frq}b$_1$-\ref{fig.frq}b$_4$.
Note that the convergences of the relative frequencies $\{\omega_i-\omega_s,~i\in\mathcal{V}_M\}$ to zero
as shown in Fig.~\ref{fig.frq}c$_1$-\ref{fig.frq}c$_4$
are the main concern of primary frequency control \cite{Dorfler2014,Motter2013}.
Since the economic power dispatch is solved on-line, it can be observed in Fig.~\ref{fig.frq}e$_1$-\ref{fig.frq}e$_2$ that
the marginal costs 
of all the controllers are the same during the transient phase under the control of PIAC and GB.
In contrast with PIAC and GB,
the marginal costs of DAI are not identical during the transient phase even 
though they achieve a consensus at the steady state. Because 
there are no control coordinations between the controllers in DecI,
the marginal costs are not identical even at the steady state as shown in Fig.~\ref{fig.frq}e$_4$. Since $u_i=\alpha_i\lambda$,
the control inputs of the PIAC method and the GB method have similar dynamics as those of their marginal costs as shown in Fig.~\ref{fig.frq}f$_1$-\ref{fig.frq}f$_2$.
Note that as shown in  Fig.~\ref{fig.frq}f$_4$, the control inputs of DecI are very close to each other because of the identical setting of $k_i$ and 
the small differences between the frequency deviations. As shown in Fig.~\ref{fig.frq}f$_1$-Fig.~\ref{fig.frq}f$_3$, the control inputs of some generators 
in the PIAC, GB and DAI control laws are small due to the high control prices.

We remark that the larger the gain coefficient $k_{GB}$ of the GB
control law the larger are the oscillations of frequency deviations even though the frequencies converge to the nominal frequency faster. 
However, the control inputs of the PIAC method converge to the power imbalance 
faster under larger control gain $k$, which leads to smaller frequency deviations. 
As shown in Fig.~\ref{fig.frq}a$_1$, the frequency drops about 0.4 Hz which can be even smaller
when $k$ is larger.
However, $k$ is related to the response time of the 
control devices and hence cannot be infinitely large. 
If the step-wise increase of the loads is too big and the gain coefficient $k$ is not large enough, i.e., the controllers 
cannot respond quickly enough, the frequency might become so low that they damage the 
synchronous machines.

\onecolumn
  
\twocolumn

\subsection{Numerical results of the Multi-area implementation}

\begin{figure}[ht]
  \begin{center}
 \includegraphics[scale=1.1]{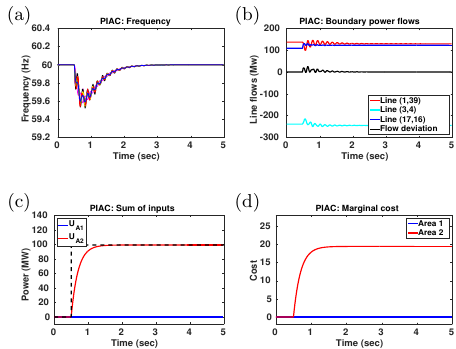}
 \caption{{Multi-area implementation of PIAC. }}
 \label{fig.multi}
 \end{center}
 \end{figure}

In this subsection, the network is separated into two areas by the red dashed line.  All the parameters of the controllers are the 
same to the ones in the single-area implementation. There are 3 generators in area $A_1$ and 7 generators in the area $A_2$ 
i.e., $\mathcal{V}_{M_1}=\{30, 37, 38\}, \mathcal{V}_{M_2}=\{31, 32, 33, 34, 35, 36, 39\}$. The boundary lines of area $A_1$ are 
in $\mathcal{E}_{A_1}=\{(1,39), (3,4), (17, 16)\}$. As in subsection \ref{subsec:numerical_single}, the secondary 
frequency controllers are installed at the nodes of the generators. 
After the step-wise increase 
of the loads at buses 4, 12 and 20 with the total amount of 99 MW in the area $A_2$, the 
multi-area implementation of PIAC recovers 
the nominal frequency as shown in Fig.~\ref{fig.multi}a and the power export deviation of area $A_1$ converges to zero as
shown by the black dashed lines in Fig.~\ref{fig.multi}b. Fig.~\ref{fig.multi}b also shows that as the system converges to a new state, the 
power flows in the three line in $\mathcal{E}_{A_1}$ are different from the ones before the step-wise increase of the loads in area $A_2$.
A characteristic of PIAC is that it decouples the control actions of the areas, which 
can be observed in Fig.~\ref{fig.multi}c. Since the step-wise increase of the power loads 
at the buses 4, 12 and 20 happen in area $A_2$, the control inputs of area $A_1$ are zero and the power is 
balanced by the controllers in the area $A_2$. This shows that with the PIAC method, the power can be balanced locally in 
an area without influencing to its neighbors. 
\emph{This characteristic of PIAC is attractive for such a non-cooperative multi-area control of a power system that 
different areas might have different amount of renewable energy. 
It is fair for the area with a large amount of renewable energy to respond to the disturbance in its own area.
As mentioned in subsection 4.2,
 this characteristic allows the controllers in different areas control the 
 system in an asynchronous way at any time according to the power imbalance within the area.}

 \section{Conclusion}\label{Sec:conclusion}

In this paper, we proposed a secondary frequency control approach, called PIAC,  to restore
the nominal frequency of power systems with a minimized control cost. A main feature  
of PIAC is that the estimated power imbalance converges exponentially with a 
speed that only depends on a gain coefficient which can be 
tuned according to the sensitivity of the control devices. Hence PIAC eliminates 
the drawback of the traditional integral control based secondary frequency control approach, 
in which large control gain coefficients lead to an overshoot problem of the control inputs and 
small ones result in a large frequency deviation of the power systems. When implemented 
in a network with multiple areas, PIAC decouples the control actions of the areas such that
the power supply and demand of an area can be balanced locally in the area without
any influences to the neighbors. For the power systems with a large amount of integrated renewable energy, 
the large transient frequency deviation can 
be reduced by PIAC with advanced control devices and communication networks. 
 
 However, in practice, there is usually some noise from the measurement of frequency and time 
 delays and even information losses in the communication. In addition, the resistance of the transmission lines
 cannot be neglected in some power networks, e.g., distribution grids of power systems or Micro-Grids. Hence, further investigation on the 
 performance of PIAC on such a lossy power network with noisy measurements and time delays
 is needed. Further investigation is also required into the performance 
 of PIAC in the power system with dynamic actuators.

 \section*{Acknowledgment}

The authors thank Prof. Chen Shen from Tsinghua University for his valuable discussions on this topic 
of the frequency control of power systems and useful advices on the simulations.  
Kaihua Xi thanks the China Scholarship Council for financial support.

\appendix

\section{Preliminaries on DAE systems}

Consider the following \emph{Differential Algebraic Equation} (DAE) systems
\begin{subequations}\label{eq:DAEtt}
\begin{align}
 \dot{x}&=f(x,y),\label{eq:DAE}\\
  0&=g(x,y), \label{Eq:DAE2}
\end{align}
\end{subequations}
where $x\in \mathbb{R}^n$ are state variables, $y\in \mathbb{R}^m$ are algebraic variables and $f:\mathbb{R}^n\times \mathbb{R}^m\rightarrow \mathbb{R}^n$ and 
$g:\mathbb{R}^n\times \mathbb{R}^m\rightarrow \mathbb{R}^m$ are 
twice continuously differentiable functions. (\ref{eq:DAE}) and (\ref{Eq:DAE2}) 
are differential and algebraic equations respectively.
$(x(x_0,y_0,t),y(x_0,y_0,t))$ is the solution with 
the admissible initial conditions $(x_0,y_0)$ satisfying the algebraic constraints
\begin{align}
 0=g(x_0,y_0). \label{eq:initial}
\end{align}
and the maximal domain of a solution of (\ref{eq:DAEtt}) is denoted by $\mathcal{I}\subset \mathbb{R}_{\geq 0}$ where 
$\mathbb{R}_{\geq 0}=\{t\in \mathbb{R}| t\geq 0\}$. 

Before presenting the Lyapunove/LaSalle stability criterion of the DAE system, we make the following two assumptions.

\textbf{Assumption A.1:} The DAE system possesses an equilibrium state $(x^*,y^*)$ such 
that $f(x^*,y^*)=0,~ g(x^*,y^*)=0$. 

\textbf{Assumption A.2:} Let $\Omega\subseteq \mathbb{R}^n\times \mathbb{R}^m$ be an open connected set containing $(x^*,y^*)$,
assume ({\ref{Eq:DAE2}}) is \emph{regular} such that 
the Jacobian of $g$ with respect to $y$ is a full rank matrix for any $(x,y)\in \Omega$, i.e.,
\begin{align*}
\text{rank}(\nabla_{y}g(x,y))=m, ~\forall (x,y)\in \Omega.
\end{align*}
Assumption A.2 ensures the existence and uniqueness of the solutions 
of (\ref{eq:DAEtt}) in $\Omega$ over the interval $\mathcal{I}$ with the 
initial condition $(x_0,y_0)$ satisfying (\ref{eq:initial}). 

The following theorem provides a sufficient stability criterion of the equilibrium of DAE in (\ref{eq:DAEtt}).

{\bf{Theorem A.1}}
\emph{(Lyapunov/LaSalle stability criterion \cite{Schiffer,Hill1990}): Consider the DAE system in (\ref{eq:DAEtt}) with 
assumptions A.1 and A.2, and an equilibrium $(x^*,y^*)\in \Omega_H\subset \Omega$. If there exists 
a continuously differentiable function $H:\Omega_H\rightarrow R$, such 
that $(x^*,y^*)$ is a strict minimum of $H$ i.e., $\nabla H|_{(x^*,y^*)}=0$ and $\nabla^2H|_{(x^*,y^*)}>0$, 
and $\dot{H}(x,y)\leq 0,~\forall (x,y)\in \Omega_H$, then the following 
statements hold:}

\emph{(1). $(x^*,y^*)$ is a stable equilibrium with a local Lyapunov function 
$V(x,y)=H(x,y)-H(x^*,y^*)\geq 0$ for $(x,y)$ near $(x^*,y^*)$,}

\emph{(2). Let $\Omega_c=\{(x,y)\in \Omega_H|H(x,y)\leq c\}$ be a compact sub-level set for 
some $c>H(x^*,y^*)$. If no solution can stay in $\{(x,y)\in \Omega_c|\dot{H}(x,y)= 0\}$ 
other than $(x^*,y^*)$, then $(x^*,y^*)$ is asymptotically stable. }

We refer to \cite{Schiffer} and \cite{Hill1990} for the proof of Theorem A.1.

  \bibliography{ifacconf.bib}         
\bibliographystyle{plain}

\end{document}